\documentclass[11pt]{article}
\usepackage{mathrsfs}
\usepackage[centertags]{amsmath}
\usepackage{amsfonts}
\usepackage{amssymb}
\usepackage{amsthm}
\usepackage{cases}
\usepackage{indentfirst}
\usepackage{epsfig}
\usepackage{graphicx}
\usepackage{color}

\date{}

\definecolor{c20}{rgb}{0.,0.7,0.}
\definecolor{c30}{rgb}{0.,0.,1.}
\definecolor{c40}{rgb}{1,0.1,0.7}
\definecolor{c50}{rgb}{1,0,0}

\def\lx#1{{\textcolor{c30}{#1}}}

\def\lx#1{#1}

\newtheorem{theorem}{Theorem}[section]

\newtheorem{lemma}{Lemma}[section]

\newtheorem{example}{Example}[section]
\numberwithin{equation}{section}

\def\P{\operatorname*{\mathbb{P}}}
\def\R{\operatorname*{\mathbb{R}}}

\def\E{\operatorname*{\mathbf{E}}}

\allowdisplaybreaks[4]
\usepackage{setspace}
\begin{document}
\title{Second-order asymptotics on distributions of maxima \\ of bivariate elliptical arrays }
\author{{  Xin Liao$^a$\quad Zhichao Weng$^b$ \quad Zuoxiang Peng$^c$\thanks{Corresponding author. Email: pzx@swu.edu.cn}}\\
 \small{$^a$Business School, University of Shanghai for Science and Technology, 200093 Shanghai, China }\\
 \small{$^b$School of Economics and Management, Fuzhou University, 350116 Fujian, China}\\
 \small{$^c$School of Mathematics and Statistics, Southwest
University, 400715 Chongqing, China}}
\maketitle

\begin{quote}
{\bf Abstract}~~Let $\{ (\xi_{ni}, \eta_{ni}), 1\leq i \leq n, n\geq 1 \}$ be a triangular array
of independent bivariate elliptical random vectors with the same distribution function as
$(S_{1}, \rho_{n}S_{1}+\sqrt{1-\rho_{n}^2}S_{2})$, $\rho_{n}\in (0,1)$, where $(S_{1},S_{2})$
is a bivariate spherical random vector. For the distribution function of radius
$\sqrt{S_{1}^2+S_{2}^2}$ belonging to the max-domain of attraction of the
Weibull distribution, Hashorva (2006) derived the limiting distribution of
maximum of this triangular array if convergence rate of $\rho_{n}$ to $1$ is given. In this paper, under the refinement of the
rate of convergence of $\rho_{n}$ to $1$ and the second-order regular variation of the distributional tail of radius,
precise second-order distributional expansions of the normalized maxima of bivariate elliptical triangular arrays are established.

{\bf Key words}~~Bivariate elliptical triangular array; Maximum;
Second-order expansion; Second-order regular variation.

{\bf AMS 2000 subject classification}~~Primary 62E20, 60G70; Secondary 60F15, 60F05.

\end{quote}

\section{Introduction}
\label{sec1}

Let
$\{(\xi_{ni},\eta_{ni}), 1\leq i \leq n, n\geq 1\}$ be a triangular array of independent bivariate elliptical random
vectors with stochastic representation
\begin{equation}\label{eq1.1}
(\xi_{ni}, \eta_{ni})\overset{d}{=} \left(S_{1},\rho_{n}S_{1}+\sqrt{1-\rho_{n}^{2}}S_{2}\right), 1\leq i \leq n, n\geq 1,
\end{equation}
where $\rho_{n}\in [-1,1]$ and $(S_1, S_2)$ be a bivariate spherical random vector with radius $R=\sqrt{S_{1}^{2}+S_{2}^{2}}$.
Here $\overset{d}{=}$ means equality in distribution.
The bivariate maxima $\mathbf{M}_{n}$ is defined componentwise by \[ \mathbf{M}_{n}=(M_{n1},M_{n2})=\left(\max_{1 \leq i\leq n} \xi_{{ni}}, \max_{1\leq i \leq n} \eta_{{ni}}\right). \]
If $S_{1}$ and $S_{2}$ are two independent standard Gaussian random variables,
\eqref{eq1.1} just represents the well-known bivariate Gaussian triangular arrays studied by
H\"{u}sler and Reiss (1989) and Kabluchko et al. (2009) who showed that
\[  \lim_{n\to \infty} \sup_{x,y\in \R} \left| \P\left( M_{n1}\leq \bar{a}_{n}x+\bar{b}_{n},
M_{n2}\leq \bar{a}_{n}y+\bar{b}_{n} \right) -K_{\lambda}(x,y) \right|=0\]
if and only if the following H\"{u}sler-Reiss condition
\[\lim_{n\to \infty} \frac{\bar{b}_{n}(1-\rho_{n})}{\bar{a}_{n}}=2\lambda^2 \in [0,\infty]\]
holds, where $K_{\lambda}(x,y)$
is the H\"{u}sler-Reiss max-stable distribution, and the norming constants
$\bar{b}_{n}$ satisfies $\sqrt{2\pi}n^{-1}\bar{b}_{n}\exp(\bar{b}_{n}^2/2)=1$ and $\bar{a}_{n}=\bar{b}_{n}^{-1}$.

Motivated by the seminal work of H\"{u}sler and Reiss (1989), Hashorva (2006) considered limit laws of maxima of bivariate elliptical triangular array with
distribution function (df) $F$ of radius $R$ belonging to the max-domain of attraction of the Weibull extreme value distribution
$\Psi_{\alpha}$, $\alpha>0$. Precisely, if the following condition
\begin{equation}\label{eq1.3}
  \lim_{n\to \infty} \frac{1-\rho_{n}}{a_{n}}=2\lambda^2 \in (0,\infty)
\end{equation}
holds with $a_{n}=1-G^{\leftarrow}\left( 1-n^{-1}¡¡\right)$, Hashorva (2006) showed that
\begin{equation}\label{eq1.2}
\lim_{n\to \infty} \sup_{(x,y)\in (-\infty,0)^2} \left| \P\{ M_{n1}\leq 1+a_{n}x, M_{n2}\leq 1+a_{n}y \}
-H_{\alpha+\frac{1}{2}, \lambda}(x,y) \right|=0,
\end{equation}
where $G$ is the df of $S_{1}$ with upper endpoint 1, $G^{\leftarrow}$ is its generalised inverse function, and
\begin{equation*}
H_{\alpha+\frac{1}{2}, \lambda}(x,y)= \exp\left(-|x|^{\alpha+\frac{1}{2}}\psi_{\alpha}
\left( \frac{1}{\sqrt{2|x|}}\left( \lambda+ \frac{y-x}{2\lambda} \right) \right)
-|y|^{\alpha+\frac{1}{2}}\psi_{\alpha}
\left( \frac{1}{\sqrt{2|y|}}\left( \lambda+ \frac{x-y}{2\lambda} \right)\right) \right).
\end{equation*}
Here, $\psi_{\alpha}$ is a df defined on $[-1,1]$ given by
\begin{equation*}
\psi_{\alpha}(z)=\frac{\Gamma(\alpha+\frac{3}{2})}{\Gamma(\alpha+1)\sqrt{\pi}}\int_{-1}^{z} (1-s^2)^{\alpha} ds,\  \alpha>0,\  z\in [-1,1],
\end{equation*}
where $\Gamma(\cdot)$ is Gamma function. For more work on extremal behaviors of extended
H\"{u}sler-Reiss bivariate Gaussian settings, see, e.g., Hooghiemstra and H\"{u}sler (1996), Hashorva (2005, 2008),
Frick and Reiss (2010), Hashorva et al. (2012), Hashorva and Weng (2013), Hashorva and Ling (2016), and Liao et al. (2016). For multivariate dependent Gaussian case, see Hashorva et al. (2015).

Recently, the studies on convergence rate of normalized maxima of H\"{u}sler-Reiss bivariate Gaussian settings and
their extensions have been received more attentions. For bivariate Gaussian triangular arrays,
Hashorva et al. (2016) established the
higher-order distributional expansions of normalized  maxima under refined  H\"{u}sler-Reiss condition, and the uniform convergence rates of normalized maxima and the second-order expansion
of the joint distributions of normalized maxima and minima were investigated by
Liao and Peng (2014, 2015), respectively. For copula version of bivariate Gaussian triangular arrays,
Frick and Reiss (2013) considered the penultimate and ultimate convergence rate for distributions of
$(n(\max_{1\leq i \leq n}\Phi(\xi_{ni})-1), n(\max_{1\leq i \leq n}\Phi(\eta_{ni})-1))$, where
$\Phi(\cdot)$ is the standard Gaussian distribution.

In this paper, we are interested in the second-order distributional expansions  of bivariate normalized
maxima as the distribution $F$ of radius $R$ belongs to the max-domain of attraction of $\Psi_{\alpha}$, $\alpha>0$. From Resnick (1987), it follows that
$F$ is in the Weibull max-domain of attraction, i.e.,
\[\lim_{n\to \infty} \sup_{x\in \R} |F^{n}(\tilde{a}_{n}x+\tilde{b}_{n})-\Psi_{\alpha}(x)|=0\]
if and only if the upper endpoint $\omega:=\sup\{ t:F(t)<1 \}<\infty$, and $1-F(\omega-s^{-1})\in RV_{-\alpha}$, regularly varying functions with exponent $-\alpha$. Furthermore, the constants can be chosen as
\[ \tilde{a}_{n}:=\omega-F^{\leftarrow}(1-n^{-1}),\  \tilde{b}_{n}:=\omega,\  n>1. \]

In order to get the desired results, we assume that $1-F(\omega-s^{-1})$ has the properties of second-order regularly variation with the
first-order parameter $-\alpha$ and the second-order parameter $\tau \leq 0$ (written as $1-F(\omega-s^{-1})\in 2RV_{-\alpha, \tau}$ ), i.e., there exists some ultimately positive or negative function $A(t)$ with $\lim_{t\to \infty} A(t)=0$ such that
\[ \lim_{t\to \infty} \frac{\frac{1-F(\omega-(ts)^{-1})}{ 1-F(\omega-t^{-1})}-s^{-\alpha}}{A(t)}
=s^{-\alpha}\frac{s^{\tau}-1}{\tau}, \quad s>0, \]
where $\frac{s^{\tau}-1}{\tau}$ is interpreted as $\log s$ when $\tau=0$, c.f., de Haan and Ferreira (2006). Furthermore, we need the rate of convergence imposed on \eqref{eq1.3} as $\lambda\in (0,\infty)$. Assume that throughout this paper that there exists $c_{n}$ with $\lim_{n\to \infty} c_{n}=0$ such that
\begin{eqnarray}\label{eq1.4}
  \lim_{n\to \infty} \frac{\lambda_{n}-\lambda}{c_{n}}=\gamma \in \R
\end{eqnarray}
holds with $\lambda_{n}=\left( \frac{1-\rho_{n}}{2a_{n}} \right)^{{1}/{2}}$ and $\lambda\in (0,\infty)$.
For two extreme cases $\lambda=0$ and $\lambda=\infty$, the analysis will be discussed
with some other additional conditions related to $\rho_{n}$.

The rest of this paper is organized as follows. In section \ref{sec2}, we provide the main results and
an example is illustrated to support our findings. Some auxiliary lemmas are given in Section \ref{sec3}, and all proofs of the
main results are deferred to Section \ref{sec4}.

\section{Main results}
\label{sec2}

In this section, we provide the main results with respect to the second-order expansions on distributions of
normalized maxima by assuming that $1-F(\omega-s^{-1})\in 2RV_{-\alpha,\tau}$, $\alpha>0, \tau\leq 0$.
There are three cases to be considered, i.e., $\lambda\in (0,\infty)$, $\lambda=0$ and $\lambda=\infty$, respectively.
For simplicity, we assume that the upper endpoint $\omega$ of $F$ is 1, which implies that the upper
endpoint of $G$ is also $1$. Throughout this paper, let $a_{n}$ be the normalized constant satisfied \lx{$a_{n}=1-G^{\leftarrow}(1-n^{-1})$}.
For the case of $\lambda\in (0,\infty)$, we need a second-order condition \eqref{eq1.4}
refining the convergence rate of condition \eqref{eq1.3}, and investigate in turn the following three cases:
$\lambda^2+x+y+\frac{(x-y)^2}{4\lambda^2}<0$, $\lambda^2+x+y+\frac{(x-y)^2}{4\lambda^2}=0$ and
$\lambda^2+x+y+\frac{(x-y)^2}{4\lambda^2}>0$, respectively. The following theorem establish the second-order
distributional expansions of  normalized maxima with $\lambda\in (0,\infty)$ and $\lambda^2+x+y+\frac{(x-y)^2}{4\lambda^2}<0$.

\begin{theorem}\label{th1}
Let $\{(\xi_{ni},\eta_{ni}), 1\leq i \leq n, n\geq 1\}$ be a triangular array satisfying \eqref{eq1.1} with $\rho_{n}\in (0,1)$. Assume that $F$ has upper endpoint equal 1 and further $1-F(1-t^{-1})\in 2RV_{-\alpha,\tau}$ with
$t>0$, $\alpha>0$, $\tau \leq 0$ and auxiliary function $A(t)$. If there exists $c_{n}$ satisfying $\lim_{n\to \infty} c_{n}=0$ such that \eqref{eq1.4}
with $\lambda_{n}=\left( \frac{1-\rho_{n}}{2a_{n}} \right)^{\frac{1}{2}}$ and $\lambda\in (0,\infty)$. With fixed $x<0$ and $y <0$ satisfying $\lambda^2 +x+y+\frac{(x-y)^2}{4\lambda^2}<0$, we have
\begin{eqnarray}\label{eq2.2}
& & \P\left( M_{n1}\leq 1+a_{n}x, M_{n2}\leq 1+a_{n}y \right)-H_{\alpha+\frac{1}{2},\lambda}(x,y) \nonumber\\
&=& H_{\alpha+\frac{1}{2},\lambda}(x,y)\Bigg\{
-\frac{1}{2}n^{-1}\left( |x|^{\alpha+\frac{1}{2}}\psi_{\alpha}\left( \frac{\lambda+\frac{y-x}{2\lambda}}{\sqrt{2|x|}} \right)
+|y|^{\alpha+\frac{1}{2}}\psi_{\alpha}
\left( \frac{\lambda+\frac{x-y}{2\lambda}}{\sqrt{2|y|}} \right) \right)^{2} \nonumber\\
& &
+Q_{\alpha+\frac{1}{2},\lambda}(x,y)
+ o(a_{n}+c_{n}+A(a_{n}^{-1})+n^{-1}) \Bigg\}
\end{eqnarray}
for large $n$, where
\begin{eqnarray*}
& & Q_{\alpha+\frac{1}{2}, \lambda}(x,y) \nonumber\\
&=& |x|^{\alpha+\frac{1}{2}}\left\{ -\frac{a_{n}|x|}{c_{\alpha}}\int_{-1}^{\frac{\lambda+\frac{y-x}{2\lambda}}{\sqrt{2|x|}}}(1-s^2)^{\alpha-1}s^{2}
(\alpha+1-\frac{3}{2}\alpha s^2-s^2)ds \right. \nonumber\\
& & \left.
-\frac{A(a_{n}^{-1})|x|^{-\tau}}{c_{\alpha}}\int_{-1}^{\frac{\lambda+\frac{y-x}{2\lambda}}{\sqrt{2|x|}}}
(1-s^2)^{\alpha}\frac{(1-s^2)^{-\tau}-1}{\tau}ds + \psi_{\alpha}\left( \frac{\lambda+\frac{y-x}{2\lambda}}{\sqrt{2|x|}} \right)
\left( - A(a_{n}^{-1})\frac{|x|^{-\tau}-1}{\tau}
\right.\right. \nonumber \\
& & \quad \left. \left.
+\frac{a_{n}}{c_{\alpha}}\int_{-1}^{1}(1-s^2)^{\alpha-1}s^2(\alpha+1-\frac{3}{2}\alpha s^2 -s^2)ds
+ \frac{A(a_{n}^{-1})}{c_{\alpha}}\int_{-1}^{1}\frac{(1-s^2)^{\alpha-\tau}-(1-s^2)^{\alpha}}{\tau} ds
\right)  \right. \nonumber \\
& & \left.
-\frac{\left( 1-\frac{\left( \lambda +\frac{y-x}{2\lambda} \right)^2}{2|x|} \right)^{\alpha}}{c_{\alpha}\sqrt{2|x|}}
\left( c_{n}\gamma\left( 1-\frac{y-x}{2\lambda^2}  \right)  -\frac{a_{n}}{2}\left( \lambda(y-x)+\frac{x(y-x)}{\lambda}
 +\frac{3(y-x)^2}{4\lambda} + \frac{(y-x)^3}{8\lambda^{3}}\right)  \right) \right\}\nonumber \\
& &
+ |y|^{\alpha+\frac{1}{2}}\left\{ -\frac{a_{n}|y|}{c_{\alpha}}\int_{-1}^{\frac{\lambda+\frac{x-y}{2\lambda}}{\sqrt{2|y|}}}(1-s^2)^{\alpha-1}s^{2}
(\alpha+1-\frac{3}{2}\alpha s^2-s^2)ds \right. \nonumber\\
& & \left.
-\frac{A(a_{n}^{-1})|y|^{-\tau}}{c_{\alpha}}\int_{-1}^{\frac{\lambda+\frac{x-y}{2\lambda}}{\sqrt{2|\lx{y}|}}}
(1-s^2)^{\alpha}\frac{(1-s^2)^{-\tau}-1}{\tau}ds + \psi_{\alpha}\left( \frac{\lambda+\frac{x-y}{2\lambda}}{\sqrt{2|y|}} \right)
\left( - A(a_{n}^{-1})\frac{|y|^{-\tau}-1}{\tau}
\right.\right. \nonumber \\
& & \quad \left. \left.
+\frac{a_{n}}{c_{\alpha}}\int_{-1}^{1}(1-s^2)^{\alpha-1}s^2(\alpha+1-\frac{3}{2}\alpha s^2 -s^2)ds
+ \frac{A(a_{n}^{-1})}{c_{\alpha}}\int_{-1}^{1}\frac{(1-s^2)^{\alpha-\tau}-(1-s^2)^{\alpha}}{\tau} ds
\right)  \right. \nonumber \\
& & \left.
-\frac{\left( 1-\frac{\left( \lambda +\frac{x-y}{2\lambda} \right)^2}{2|y|} \right)^{\alpha}}{c_{\alpha}\sqrt{2|y|}}
\left( c_{n}\gamma\left( 1-\frac{x-y}{2\lambda^2}  \right)  +\frac{a_{n}}{2}\left( \lambda(y-x)+\frac{x(y-x)}{\lambda}
 +\frac{(y-x)^2}{4\lambda} + \frac{(y-x)^3}{8\lambda^{3}}\right)  \right) \right\} \nonumber \\
& &
+o\left( a_{n}+c_{n}+A(a_{n}^{-1}) \right),
\end{eqnarray*}
where $c_{\alpha}=\int_{0}^{1}(1-s)^{\alpha}s^{-\frac{1}{2}}ds=\frac{\Gamma(\alpha+1)\sqrt{\pi}}
{\Gamma(\alpha+\frac{3}{2})}$.
\end{theorem}

For the case of $\lambda\in (0,\infty)$  with $\lambda^2+x+y+\frac{(x-y)^2}{4\lambda^2}=0$, restrictions on $x$ and $y$ may be needed.  There are three cases to be considered: i), \lx{$x<0$ with $y=-\left( \sqrt{2}\lambda-\sqrt{-x} \right)^2$, and $y<0$ with $x=-\left( \sqrt{2}\lambda - \sqrt{-y} \right)^2$}; ii), $x<0$ \lx{with} $y=-(\sqrt{2}\lambda+\sqrt{-x})^{2}$ and, iii), $y<0$ \lx{with} $x=-(\sqrt{2}\lambda+\sqrt{-y})^{2}$. The asymptotics are different and we present the results in turn.

\begin{theorem}\label{th2}
Let $\{(\xi_{ni},\eta_{ni}), 1\leq i \leq n, n\geq 1\}$ be a triangular array satisfying \eqref{eq1.1} with $\rho_{n}\in (0,1)$. Assume that $F$ has upper endpoint equal 1 and further $1-F(1-t^{-1})\in 2RV_{-\alpha,\tau}$ with $t>0$,
$\alpha>0$, $\tau \leq 0$ and auxiliary function $A(t)$. Assuming that \eqref{eq1.4} holds with $\lambda\in (0,\infty)$. \lx{For $x<0$ with $y=-\left( \sqrt{2}\lambda-\sqrt{-x} \right)^2$, and $y<0$ with $x=-\left( \sqrt{2}\lambda - \sqrt{-y} \right)^2$}, we have
\begin{itemize}
\item[(i)]~~if $\lim_{n\to \infty} \frac{c_{n}}{a_{n}}=k \in \R$, then
\begin{eqnarray}\label{eq2.3}
& & \P\left(M_{n1}\leq 1+a_{n}x,M_{n2}\leq 1+a_{n}y\right) -H_{\alpha+\frac{1}{2},\lambda}(x,y) \nonumber\\
&=& H_{\alpha+\frac{1}{2},\lambda}(x,y) \left\{ -|x|^{\alpha+\frac{1}{2}}\left( \frac{a_{n}(|x|-1)}{c_{\alpha}}
\int_{0}^{1}(1-s)^{\alpha-1}s^{\frac{1}{2}}\left( \alpha+1-\frac{3}{2}\alpha s -s \right)ds \right. \right. \nonumber\\
& & \quad \left.  \left.
+A(a_{n}^{-1})\left( \frac{|x|^{-\tau}-1}{\tau}+\frac{|x|^{-\tau}-1}{c_{\alpha}}
\int_{0}^{1}\frac{(1-s)^{\alpha-\tau}-(1-s)^{\alpha}}{\tau}s^{-\frac{1}{2}}ds  \right)  \right) \right. \nonumber\\
& &
\left.  -|y|^{\alpha+\frac{1}{2}}\left( \frac{a_{n}(|y|-1)}{c_{\alpha}}
\int_{0}^{1}(1-s)^{\alpha-1}s^{\frac{1}{2}}\left( \alpha+1-\frac{3}{2}\alpha s -s \right)ds \right. \right. \nonumber\\
& & \quad \left.  \left.
+A(a_{n}^{-1})\left( \frac{|y|^{-\tau}-1}{\tau}+\frac{|y|^{-\tau}-1}{c_{\alpha}}
\int_{0}^{1}\frac{(1-s)^{\alpha-\tau}-(1-s)^{\alpha}}{\tau}s^{-\frac{1}{2}}ds  \right)  \right) \right. \nonumber \\
& &
\left. -\frac{1}{2}n^{-1}\left( |x|^{\alpha+\frac{1}{2}} + |y|^{\alpha+\frac{1}{2}} \right)^{2}
+ o\left( a_{n}+n^{-1}+A(a_{n}^{-1}) \right) \right\}
\end{eqnarray}
for large $n$;

\item[(ii)]~~if $a_{n}=o(c_{n})$, then
\begin{eqnarray}\label{eq2.4}
& & \P\left( M_{n1}\leq 1+a_{n}x, M_{n2}\leq 1+a_{n}y \right) - H_{\alpha+\frac{1}{2}, \lambda}(x,y)\nonumber\\
&=&H_{\alpha+\frac{1}{2}, \lambda}(x,y)\left\{-|x|^{\alpha+\frac{1}{2}}A(a_{n}^{-1})\left(
\frac{|x|^{-\tau}-1}{\tau} +\frac{|x|^{-\tau}-1}{c_{\alpha}}\int_{0}^{1}\frac{(1-s)^{\alpha-\tau}-(1-s)^{\alpha}}{\tau}s^{-\frac{1}{2}}ds \right) \right. \nonumber\\
& &
\left. - |y|^{\alpha+\frac{1}{2}}A(a_{n}^{-1})\left(
\frac{|y|^{-\tau}-1}{\tau} +\frac{|y|^{-\tau}-1}{c_{\alpha}}\int_{0}^{1}\frac{(1-s)^{\alpha-\tau}-(1-s)^{\alpha}}{\tau}s^{-\frac{1}{2}}ds \right)  \right. \nonumber \\
& &
\left. -\frac{1}{2}n^{-1}\left( |x|^{\alpha+\frac{1}{2}} + |y|^{\alpha+\frac{1}{2}} \right)^{2}
+ o\left( c_{n}+n^{-1}+A(a_{n}^{-1}) \right) \right\}
\end{eqnarray}
for large $n$.

\end{itemize}
\end{theorem}

\begin{theorem}\label{th3}
Let $\{(\xi_{ni},\eta_{ni}), 1\leq i \leq n, n\geq 1\}$ be a triangular array satisfying \eqref{eq1.1} with $\rho_{n}\in (0,1)$. Assume that $F$ has upper endpoint equal 1 and further $1-F(1-t^{-1})\in 2RV_{-\alpha,\tau}$ with $t>0$
$\alpha>0$, $\tau \leq 0$ and auxiliary function $A(t)$. Assuming that \eqref{eq1.4} holds with $\lambda\in (0,\infty)$. \lx{For} $x<0$ \lx{with} $y=-(\sqrt{2}\lambda+\sqrt{-x})^{2}$, we have
\begin{itemize}
\item[(i)]~~if $\lim_{n\to \infty}\frac{c_{n}}{a_{n}}=k \in \R$, then
\begin{eqnarray}\label{eq2.5}
& & \P\left(M_{n1}\leq 1+a_{n}x,M_{n2}\leq 1+a_{n}y\right) -H_{\alpha+\frac{1}{2},\lambda}(x,y) \nonumber\\
&=& H_{\alpha+\frac{1}{2},\lambda}(x,y) \left\{ -|y|^{\alpha+\frac{1}{2}}\left( \frac{a_{n}(|y|-1)}{c_{\alpha}}
\int_{0}^{1}(1-s)^{\alpha-1}s^{\frac{1}{2}}\left( \alpha+1-\frac{3}{2}\alpha s -s \right)ds \right. \right. \nonumber\\
& & \quad \left.  \left.
+A(a_{n}^{-1})\left( \frac{|y|^{-\tau}-1}{\tau}+\frac{|y|^{-\tau}-1}{c_{\alpha}}
\int_{0}^{1}\frac{(1-s)^{\alpha-\tau}-(1-s)^{\alpha}}{\tau}s^{-\frac{1}{2}}ds  \right)  \right) \right. \nonumber \\
& &
\left. -\frac{1}{2}n^{-1} |y|^{2\alpha+1}
+ o\left( a_{n}+n^{-1}+A(a_{n}^{-1}) \right) \right\}
\end{eqnarray}
for large $n$;

\item[(ii)]~~if $a_{n}=o(c_{n})$, then
\begin{eqnarray}\label{eq2.6}
& & \P\left(M_{n1}\leq 1+a_{n}x,M_{n2}\leq 1+a_{n}y\right) -H_{\alpha+\frac{1}{2},\lambda}(x,y) \nonumber\\
&=& H_{\alpha+\frac{1}{2},\lambda}(x,y) \left\{-|y|^{\alpha+\frac{1}{2}}
A(a_{n}^{-1})\left( \frac{|y|^{-\tau}-1}{\tau}+\frac{|y|^{-\tau}-1}{c_{\alpha}}
\int_{0}^{1}\frac{(1-s)^{\alpha-\tau}-(1-s)^{\alpha}}{\tau}s^{-\frac{1}{2}}ds  \right) \right. \nonumber\\
& & \left.
-\frac{1}{2}n^{-1} |y|^{2\alpha+1}
+ o\left( c_{n}+n^{-1}+A(a_{n}^{-1}) \right) \right\}
\end{eqnarray}
for large $n$.

\end{itemize}
\end{theorem}

\begin{theorem}\label{th4}
Let $\{(\xi_{ni},\eta_{ni}), 1\leq i \leq n, n\geq 1\}$ be a triangular array satisfying \eqref{eq1.1} with $\rho_{n}\in (0,1)$. Assume that $F$ has upper endpoint equal 1 and further $1-F(1-t^{-1})\in 2RV_{-\alpha,\tau}$ with $t>0$
$\alpha>0$, $\tau \leq 0$ and auxiliary function $A(t)$. Assuming that \eqref{eq1.4} holds with $\lambda\in (0,\infty)$. \lx{For} $y<0$ \lx{with} $x=-(\sqrt{2}\lambda+\sqrt{-y})^{2}$, we have
\begin{itemize}
\item[(i)]~~if $\lim_{n\to \infty}\frac{c_{n}}{a_{n}}=k \in \R$, then
\begin{eqnarray}\label{eq2.7}
& & \P\left(M_{n1}\leq 1+a_{n}x,M_{n2}\leq 1+a_{n}y\right) -H_{\alpha+\frac{1}{2},\lambda}(x,y) \nonumber\\
&=& H_{\alpha+\frac{1}{2},\lambda}(x,y) \left\{-|x|^{\alpha+\frac{1}{2}}\left( \frac{a_{n}(|x|-1)}{c_{\alpha}}
\int_{0}^{1}(1-s)^{\alpha-1}s^{\frac{1}{2}}\left( \alpha+1-\frac{3}{2}\alpha s -s \right)ds \right. \right. \nonumber\\
& & \quad \left.  \left.
+A(a_{n}^{-1})\left( \frac{|x|^{-\tau}-1}{\tau}+\frac{|x|^{-\tau}-1}{c_{\alpha}}
\int_{0}^{1}\frac{(1-s)^{\alpha-\tau}-(1-s)^{\alpha}}{\tau}s^{-\frac{1}{2}}ds  \right)  \right) \right. \nonumber \\
& &
\left. -\frac{1}{2}n^{-1} |x|^{2\alpha+1}
+ o\left( a_{n}+n^{-1}+A(a_{n}^{-1}) \right) \right\}
\end{eqnarray}
for large $n$;

\item[(ii)]~~if $a_{n}=o(c_{n})$, then
\begin{eqnarray}\label{eq2.8}
& & \P\left(M_{n1}\leq 1+a_{n}x,M_{n2}\leq 1+a_{n}y\right) -H_{\alpha+\frac{1}{2},\lambda}(x,y) \nonumber\\
&=& H_{\alpha+\frac{1}{2},\lambda}(x,y) \left\{ -|x|^{\alpha+\frac{1}{2}}
A(a_{n}^{-1})\left( \frac{|x|^{-\tau}-1}{\tau}+\frac{|x|^{-\tau}-1}{c_{\alpha}}
\int_{0}^{1}\frac{(1-s)^{\alpha-\tau}-(1-s)^{\alpha}}{\tau}s^{-\frac{1}{2}}ds  \right) \right. \nonumber\\
& & \left.
-\frac{1}{2}n^{-1} |x|^{2\alpha+1}
+ o\left( c_{n}+n^{-1}+A(a_{n}^{-1}) \right) \right\}
\end{eqnarray}
for large $n$.

\end{itemize}
\end{theorem}

For $x<0,y<0$, note that $\lambda^2+x+y+\frac{(x-y)^2}{4\lambda^2}>0$ ensures that $\frac{\left( \lambda+\frac{y-x}{2\lambda} \right)^2}{-2x}>1$ and $\frac{\left( \lambda+\frac{x-y}{2\lambda} \right)^2}{-2y}>1$ all hold. Since
$\lambda+\frac{y-x}{2\lambda}<-\sqrt{-2x}$ and $\lambda+\frac{x-y}{2\lambda}<-\sqrt{-2y}$ have no intersection,
we can deal with the case of $\lambda \in (0,\infty)$ and $\lambda^2+x+y+\frac{(x-y)^2}{4\lambda^2}>0$ by the following three cases:
(i) $\lambda+\frac{y-x}{2\lambda}>\sqrt{-2x}, \lambda+\frac{x-y}{2\lambda}>\sqrt{-2y}$;
(ii) $\lambda+\frac{y-x}{2\lambda}<-\sqrt{-2x}, \lambda+\frac{x-y}{2\lambda}>\sqrt{-2y}$;
(iii) $\lambda+\frac{y-x}{2\lambda}>\sqrt{-2x}, \lambda+\frac{x-y}{2\lambda}< -\sqrt{-2y}$.

\begin{theorem}\label{th5}
Let $\{(\xi_{ni},\eta_{ni}), 1\leq i \leq n, n\geq 1\}$ be a triangular array satisfying \eqref{eq1.1} with $\rho_{n}\in (0,1)$. Assume that $F$ has upper endpoint equal 1 and further $1-F(1-t^{-1})\in 2RV_{-\alpha,\tau}$ with $t>0$,
$\alpha>0$, $\tau \leq 0$ and auxiliary function $A(t)$. If \eqref{eq1.3} holds with $\lambda\in (0,\infty)$, then
\begin{itemize}
\item[(i)]~~if $\lambda+\frac{y-x}{2\lambda}>\sqrt{-2x}, \lambda+\frac{x-y}{2\lambda}>\sqrt{-2y}$,
then \eqref{eq2.3} holds;

\item[(ii)]~~if $\lambda+\frac{y-x}{2\lambda}<-\sqrt{-2x}, \lambda+\frac{x-y}{2\lambda}>\sqrt{-2y}$,
then \eqref{eq2.5} holds;

\item[(iii)]~~if $\lambda+\frac{y-x}{2\lambda}>\sqrt{-2x}, \lambda+\frac{x-y}{2\lambda}< -\sqrt{-2y}$,
then \eqref{eq2.7} holds.

\end{itemize}
\end{theorem}

To end this section, under some other additional condition related to $\rho_{n}$, we can investigate two extreme cases: $\lambda=0$ and $\lambda=\infty$. For the case of $\lambda=0$, we have
the following results.

\begin{theorem}\label{th6}
Let $\{(\xi_{ni},\eta_{ni}), 1\leq i \leq n, n\geq 1\}$ be a triangular array satisfying \eqref{eq1.1} with $\rho_{n}\in (0,1]$. Assume that $F$ has upper endpoint equal 1 and further $1-F(1-t^{-1})\in 2RV_{-\alpha,\tau}$ with $t>0$,
$\alpha>0$, $\tau \leq 0$ and auxiliary function $A(t)$. Then,
\begin{itemize}
\item[(i),]~~if $\rho_{n}=1$, then{\small
\begin{eqnarray}\label{eq2.9}
& & \P\left(M_{n1}\leq 1+a_{n}x,M_{n2}\leq 1+a_{n}y\right) -H_{\alpha+\frac{1}{2},0}(x,y) \nonumber\\
&=& H_{\alpha+\frac{1}{2},0}(x,y) \left\{ -|\min(x,y)|^{\alpha+\frac{1}{2}}\left( \frac{a_{n}(|\min(x,y)|-1)}{c_{\alpha}}
\int_{0}^{1}(1-s)^{\alpha-1}s^{\frac{1}{2}}\left( \alpha+1-\frac{3}{2}\alpha s -s \right)ds \right. \right. \nonumber\\
& &  \left.  \left.
+A(a_{n}^{-1})\left( \frac{|\min(x,y)|^{-\tau}-1}{\tau}+\frac{|\min(x,y)|^{-\tau}-1}{c_{\alpha}}
\int_{0}^{1}\frac{(1-s)^{\alpha-\tau}-(1-s)^{\alpha}}{\tau}s^{-\frac{1}{2}}ds  \right)  \right) \right. \nonumber \\
& &
\left. -\frac{1}{2}n^{-1} |\min(x,y)|^{2\alpha+1}
+ o\left( a_{n}+n^{-1}+A(a_{n}^{-1}) \right) \right\}
\end{eqnarray}}
for large $n$;

\item[(ii),]~~if $\rho_{n}\in (0,1)$, further assume that $\lim_{n\to \infty} \frac{1-\rho_{n}}{a_{n}}=0$ and
$\lim_{n\to \infty}\frac{a_{n}^{{5}/{3}}}{1-\rho_{n}}=0$,
then \eqref{eq2.9} also holds.

\end{itemize}
\end{theorem}

For the case of $\lambda=\infty$, we have the following results.

\begin{theorem}\label{th7}
Let $\{(\xi_{ni},\eta_{ni}), 1\leq i \leq n, n\geq 1\}$ be a triangular array satisfying \eqref{eq1.1} with $\rho_{n}\in [-1,1)$. Assume that $F$ has upper endpoint equal 1 and further $1-F(1-t^{-1})\in 2RV_{-\alpha,\tau}$ with $t>0$,
$\alpha>0$, $\tau \leq 0$. Then for $\lambda =\infty$,
\begin{itemize}
\item[(i)]~~if $\rho_{n}\to 1$, further assume that $\lim_{n\to \infty}\frac{1-\rho_{n}}{a_{n}}=\infty$ and
$\lim_{n\to \infty} \frac{(1-\rho_n)^3}{a_{n}}=0$, then \eqref{eq2.3} holds;

\item[(ii)]~~if $\rho_{n}\to d \in[-1,1)$,
then \eqref{eq2.3} also holds.

\end{itemize}
\end{theorem}

\begin{example}
Consider a triangular array  $\{(\xi_{ni},\eta_{ni}), 1\leq i \leq n, n\geq 1\}$
satisfying \eqref{eq1.1} with almost surely positive random radius $R$ being
Beta distributed with parameters $a,b>0$. One can check that
\[ 1-F(1-x^{-1})=\bar{F}(1-x^{-1})=\frac{1}{bB(a,b)}x^{-b}\left( 1-\frac{b(a-1)x^{-1}}{b+1}(1+o(1)) \right),\ x\to \infty. \]
Hence, $\bar{F}\in 2RV_{-b,-1}$ with auxiliary function $A(t)=\frac{b(a-1)}{b+1}t^{-1}$.
Note that $a_{n}=1-G^{\leftarrow}(1-n^{-1})$. By using Lemma \ref{lem2}, we have
\begin{eqnarray*}
n^{-1}&=& 1-G(1-a_{n})\\
&=& \frac{a_{n}^{b+\frac{1}{2}}c_{b}}{\sqrt{2}\pi b B(a,b)} \left( 1+\left( \frac{1}{c_{b}}\int_{0}^{1}(1-s)^{b-1}s^{\frac{1}{2}}(b+1-\frac{3}{2}bs-s)ds
-\frac{b(a-1)}{b+1} \right. \right. \\
& & \left. \left.+\frac{b(a-1)}{(b+1)c_{b}}\int_{0}^{1}(1-s)^bs^{\frac{1}{2}}ds  \right)a_{n}
+o(a_{n})  \right)
\end{eqnarray*}
for large $n$, where $c_{b}=\int_{0}^{1} (1-s)^{b}s^{-\frac{1}{2}}ds$. So we can get
\[a_{n}=\left(  \frac{bB(a,b)\sqrt{2}\pi}{nc_{b}} \right)^{\frac{1}{b+\frac{1}{2}}}
-\frac{m}{b+\frac{1}{2}}\left( \frac{b B(a,b) \sqrt{2} \pi}{nc_{b}} \right)^{\frac{2}{b+\frac{1}{2}}}(1+o(1))\]
with
\[\lx{m= \frac{1}{c_{b}}\int_{0}^{1}(1-s)^{b-1}s^{\frac{1}{2}}(b+1-\frac{3}{2}s-s)ds
-\frac{b(a-1)}{b+1} + \frac{b(a-1)}{(b+1)c_{b}}\int_{0}^{1} (1-s)^{b} s^{\frac{1}{2}}ds .}\]
Assume that $\rho_{n}=1-n^{-\frac{1}{b+\frac{1}{2}}}$, then the second-order condition \eqref{eq1.4} holds
with\[\lambda=\frac{1}{\sqrt{2}}\left( \frac{b B(a,b) \sqrt{2} \pi}{c_{b}} \right)^{-\frac{1}{2b+1}},\quad
c_{n}=n^{-\frac{1}{b+\frac{1}{2}}},\quad \mbox{and}\quad\gamma=\frac{m}{2\sqrt{2}(b+\frac{1}{2})}\left( \frac{b B(a,b) \sqrt{2} \pi}{c_{b}} \right)^{\frac{1}{2b+1}}.\]
\lx{For $x<0,y<0$ with $\lambda^2+x+y+\frac{(x-y)^2}{4\lambda^2}\leq 0$, from Theorem \ref{th1}-\ref{th4},}  it follows that
\begin{eqnarray*}
& & \P\left( M_{n1}\leq 1+a_{n}x, M_{n2}\leq 1+a_{n}y \right) -H_{b,\lambda}(x,y)\\
&=& H_{b,\lambda}(x,y)
\left\{ \left( \frac{bB(a,b)\sqrt{2}\pi}{nc_{b}} \right)^{\frac{1}{b+\frac{1}{2}}}\left[
|x|^{b+\frac{1}{2}}\left( -\frac{|x|}{c_{b}}\int_{-1}^{\frac{\lambda+\frac{y-x}{2\lambda}}{\sqrt{2|x|}}}
(1-s^2)^{b-1}s^{2}\left( b+1-\frac{3}{2}bs^2-s^2 \right)ds \right. \right. \right. \\
& & \left. \left. \left. +\frac{b(a-1)}{c_{b}(b+1)}|x| \int_{-1}^{\frac{\lambda+\frac{y-x}{2\lambda}}{\sqrt{2|x|}}} \left( (1-s^2)^{b+1} -(1-s^2)^{b} \right)ds
+\psi_{b}\left( \frac{\lambda+\frac{y-x}{2\lambda}}{\sqrt{2|x|}}\right)\left( \frac{b(a-1)}{b+1}(|x|-1)  \right. \right. \right. \right. \\
& &  \left.  \left.  \left.  \left.   +\frac{1}{c_{b}}\int_{\lx{-1}}^{1}(1-s^2)^{b-1}s^2\left(b+1-\frac{3}{2}bs^2-s^2  \right)ds
-\frac{b(a-1)}{c_{b}(b+1)}\int_{\lx{-1}}^{1}\left( (1-s^2)^{b+1}-(1-s^2)^b \right)ds \right)  \right. \right. \right. \\
& & \left. \left. \left. -\frac{1}{c_{b}\sqrt{2|x|}}\left( 1-\frac{\left( \lambda+\frac{y-x}{2\lambda} \right)^2}{2|x|} \right)^{b}\left( \gamma\left( 1-\frac{y-x}{2\lx{\lambda^2}} \right)\left(  \frac{c_{b}}{bB(a,b)\sqrt{2}\pi} \right)^{\frac{1}{b+\frac{1}{2}}}\right.\right.\right.\right.\\
& &  \left. \left. \left. \left.
-\frac{1}{2}\left( \lambda(y-x)+\frac{(y-x)^3}{8\lambda^3} +\frac{(y-x)^2}{4\lambda} +\frac{x(y-x)}{\lambda}\right)  \right)   \right)  \right. \right. \\
& & \left. \left. + |y|^{b+\frac{1}{2}}\left( -\frac{|y|}{c_{b}}\int_{-1}^{\frac{\lambda+\frac{x-y}{2\lambda}}{\sqrt{2|y|}}}
(1-s^2)^{b-1}s^{2}\left( b+1-\frac{3}{2}bs^2-s^2 \right)ds
  \right. \right.  \right. \\
& & \left. \left.\left. +\frac{b(a-1)}{c_{b}(b+1)}|y| \int_{-1}^{\frac{\lambda+\frac{x-y}{2\lambda}}{\sqrt{2|y|}}} \left( (1-s^2)^{b+1} -(1-s^2)^{b} \right)ds
+\psi_{b}\left( \frac{\lambda+\frac{x-y}{2\lambda}}{\sqrt{2|y|}}\right)\left( \frac{b(a-1)}{b+1}(|y|-1)  \right.\right.  \right. \right. \\
& &  \left. \left. \left. \left. +\frac{1}{c_{b}}\int_{\lx{-1}}^{1}(1-s^2)^{b-1}s^2\left(b+1-\frac{3}{2}bs^2-s^2  \right)ds
-\frac{b(a-1)}{c_{b}(b+1)}\int_{\lx{-1}}^{1}\left( (1-s^2)^{b+1}-(1-s^2)^b \right)ds \right) \right. \right. \right.\\
& & \left. \left. \left. -\frac{1}{c_{b}\sqrt{2|y|}}\left( 1-\frac{\left( \lambda+\frac{x-y}{2\lambda} \right)^2}{2|y|} \right)^{b}\left( \gamma\left( 1-\lx{\frac{x-y}{2\lambda^2}} \right)\left(  \frac{c_{b}}{bB(a,b)\sqrt{2}\pi} \right)^{\frac{1}{b+\frac{1}{2}}}\right.\right.\right.\right.\\
& &  \left. \left. \left. \left.
\lx{+} \frac{1}{2}\left( \lambda(y-x)+\frac{(y-x)^3}{8\lambda^3} +\frac{(y-x)^2}{4\lambda} +\frac{x(y-x)}{\lambda}\right)  \right)   \right) \right] \right. \\
& & \left. -\frac{1}{2}n^{-1} \left( |x|^{b+\frac{1}{2}}\psi_{b}\left( \frac{\lambda+\frac{y-x}{2\lambda}}{\sqrt{2|x|}} \right)
+ |y|^{b+\frac{1}{2}}\psi_{b}\left( \frac{\lambda+\frac{x-y}{2\lambda}}{\sqrt{2|y|}}\right) \right)^2 + o\left( n^{-1}+n^{-\frac{1}{b+\frac{1}{2}}} \right) \right\}.
\end{eqnarray*}
\lx{
For $x<0$, $y<0$ with $\lambda^2+x+y+\frac{(x-y)^2}{4\lambda^2}>0$, using Theorem \ref{th5} we can get
\begin{eqnarray*}
& & \P\left( M_{n1}\leq 1+a_{n}x, M_{n2}\leq 1+a_{n}y \right) -H_{b,\lambda}(x,y)\\
&=& H_{b,\lambda}(x,y)
\left\{ \left( \frac{bB(a,b)\sqrt{2}\pi}{nc_{b}} \right)^{\frac{1}{b+\frac{1}{2}}}\left[
|x|^{b+\frac{1}{2}}\left( -\frac{|x|}{c_{b}}\int_{-1}^{sgn\left(\frac{\lambda+\frac{y-x}{2\lambda}}{\sqrt{2|x|}}\right)}
(1-s^2)^{b-1}s^{2}\left( b+1-\frac{3}{2}bs^2-s^2 \right)ds \right. \right. \right. \\
& & \left. \left. \left. +\frac{b(a-1)}{c_{b}(b+1)}|x| \int_{-1}^{sgn\left(\frac{\lambda+\frac{y-x}{2\lambda}}{\sqrt{2|x|}}\right)} \left( (1-s^2)^{b+1} -(1-s^2)^{b} \right)ds
+\psi_{b}\left( \frac{\lambda+\frac{y-x}{2\lambda}}{\sqrt{2|x|}}\right)\left( \frac{b(a-1)}{b+1}(|x|-1)  \right. \right. \right. \right. \\
& &  \left.  \left.  \left.  \left.   +\frac{1}{c_{b}}\int_{-1}^{1}(1-s^2)^{b-1}s^2\left(b+1-\frac{3}{2}bs^2-s^2  \right)ds
-\frac{b(a-1)}{c_{b}(b+1)}\int_{-1}^{1}\left( (1-s^2)^{b+1}-(1-s^2)^b \right)ds \right)    \right)  \right. \right. \\
& & \left. \left. + |y|^{b+\frac{1}{2}}\left( -\frac{|y|}{c_{b}}\int_{-1}^{sgn\left(\frac{\lambda+\frac{x-y}{2\lambda}}{\sqrt{2|y|}}\right)}
(1-s^2)^{b-1}s^{2}\left( b+1-\frac{3}{2}bs^2-s^2 \right)ds
  \right. \right.  \right. \\
& & \left. \left.\left. +\frac{b(a-1)}{c_{b}(b+1)}|y| \int_{-1}^{sgn\left(\frac{\lambda+\frac{x-y}{2\lambda}}{\sqrt{2|y|}}\right)} \left( (1-s^2)^{b+1} -(1-s^2)^{b} \right)ds
+\psi_{b}\left( \frac{\lambda+\frac{x-y}{2\lambda}}{\sqrt{2|y|}}\right)\left( \frac{b(a-1)}{b+1}(|y|-1)  \right.\right.  \right. \right. \\
& &  \left. \left. \left. \left. +\frac{1}{c_{b}}\int_{-1}^{1}(1-s^2)^{b-1}s^2\left(b+1-\frac{3}{2}bs^2-s^2  \right)ds
-\frac{b(a-1)}{c_{b}(b+1)}\int_{-1}^{1}\left( (1-s^2)^{b+1}-(1-s^2)^b \right)ds \right) \right)  \right] \right. \\
& & \left. -\frac{1}{2}n^{-1} \left( |x|^{b+\frac{1}{2}}\psi_{b}\left( \frac{\lambda+\frac{y-x}{2\lambda}}{\sqrt{2|x|}} \right)
+ |y|^{b+\frac{1}{2}}\psi_{b}\left( \frac{\lambda+\frac{x-y}{2\lambda}}{\sqrt{2|y|}}\right) \right)^2 + o\left( n^{-1}+n^{-\frac{1}{b+\frac{1}{2}}} \right) \right\}
\end{eqnarray*}
with the condition \eqref{eq1.3} holding, where sgn($\cdot$) is sign function. }
\end{example}

\section{Auxiliary lemmas}
\label{sec3}

In order to prove our main results, we give some auxiliary lemmas in this section. For simplicity, in the sequel,
let $t_{n}(x)=1+a_{n}x$, where $x<0$ and $a_{n}=1-G^{-1}(1-n^{-1})$. The first lemma is about Drees' type
inequalities for the second-order regular varying functions $2RV_{\alpha, \tau}$, cf. de Haan and Ferreira (2006).

\begin{lemma}\label{lem1}
If $\chi \in 2RV_{\alpha, \tau}$ with auxiliary function $A_{1}(t), \alpha\in R$ and $\tau \leq 0$, then for
any $\varepsilon, \delta >0$, there exist an auxiliary function $A(t), A(t)\sim A_{1}(t)$ as $t\to \infty$,
and $t_{0}=t_{0}(\varepsilon, \delta)>0$ such that for all $t, tx>t_{0}$,
\begin{eqnarray*}
\left| \frac{\frac{\chi(tx)}{\chi(t)}-x^{\alpha}}{A(t)}-x^{\alpha}\frac{x^{\tau}-1}{\tau} \right|
\leq \varepsilon x^{\alpha+\tau}\max\left( x^{\delta}, x^{-\delta} \right).
\end{eqnarray*}
\end{lemma}

Without loss of generality, we assume that auxiliary functions of $2RV$ functions are positive eventually
in the following proofs. Before providing the distribution tail expansion of $G$, we need the following lemma.

\begin{lemma}\label{lem2}
Let $Q$ be a Beta distributed random variable with positive parameters $a$ and $b$. Assume that distribution
function $F$ with upper endpoint 1 satisfies $1-F(1-t^{-1})\in 2RV_{-\alpha, \tau}$ with $\alpha>0, \tau\leq 0$
and auxiliary function $A(t)$.
Then for large $n$ we have
\begin{eqnarray}\label{eq3.1}\small
& &
\E\left( 1-F\left(t_{n}(x)(1-Q)^{-\frac{1}{2}}\right) \right)\nonumber\\
&=&\frac{(2a_{n}|x|)^a\Gamma(a+b)}
{\Gamma(a)\Gamma(b)}[1-F(t_{n}(x))]\left( \int_{0}^{1}(1-s)^{\alpha}s^{a-1}ds
+A(a_{n}^{-1})|x|^{-\tau}\int_{0}^{1}(1-s)^{\alpha}\frac{(1-s)^{-\tau}-1}{\tau}s^{a-1}ds
\right. \nonumber\\
& & \left.
+a_{n}|x|\int_{0}^{1}(1-s)^{\alpha-1}s^{a}\left( \alpha -2(b-1)-\frac{3}{2}\alpha s + 2(b-1)s \right)ds
+o\left( a_{n}+A(a_{n}^{-1}) \right)
 \right)
\end{eqnarray}
with $x<0$. Here, $\frac{(1-y)^{-\tau}-1}{\tau}$ is interpreted as $-\log (1-y)$ when $\tau=0$.
\end{lemma}

\noindent
{\bf Proof.}~In order to get the desired result, we first give the following inequalities for $0<x<\frac{1}{2}$:
\begin{eqnarray}\label{eq3.2}
1-\alpha x < (1-x)^{\alpha} \leq 1-\alpha x+\frac{\alpha (\alpha-1)}{2} x^{2},\quad \alpha \geq 2,
\end{eqnarray}
\begin{eqnarray}\label{eq3.3}
1-\alpha x \leq (1-x)^{\alpha} \lx{\leq}1-\alpha x +2^{1-\alpha}\alpha (\alpha-1) x^{2}, \quad 1\leq \alpha <2,
\end{eqnarray}
\begin{eqnarray} \label{eq3.4}
1-\alpha x +2^{1-\alpha}\alpha (\alpha-1)x^{2}< (1-x)^{\alpha} <1-\alpha x,\quad  0<\alpha<1
\end{eqnarray}
and
\begin{eqnarray}\label{eq3.5}
1+\alpha x < (1+x)^{\alpha} \leq 1+ \alpha x+\left(\frac{3}{2}\right)^{\alpha-2}\frac{\alpha (\alpha-1)}{2} x^{2},\quad \alpha \geq 2,
\end{eqnarray}
\begin{eqnarray}\label{eq3.6}
1+\alpha x \leq (1+x)^{\alpha} \lx{\leq} 1+ \alpha x+\frac{\alpha (\alpha-1)}{2} x^{2},\quad 1\leq \alpha < 2,
\end{eqnarray}
\begin{eqnarray}\label{eq3.7}
1+ \alpha x+\frac{\alpha (\alpha-1)}{2} x^{2}< (1+x)^{\alpha} < 1+\alpha x, \quad 0< \alpha < 1.
\end{eqnarray}

Note that for $\beta >-1, 1<c<\alpha+1$ and large $n$, we have
\begin{eqnarray}\label{eq3.8}
& & \int_{1-\frac{(a_{n}|x|)^{\frac{1}{c}}}{2}}^{1-\frac{a_{n}|x|}{2}}
\left( 1-t_{n}(x)s-\frac{3}{2}a_{n}|x|t_{n}(x)s^{2}  \right)^{\alpha} s^{\beta}ds \nonumber \\
&=& \int_{1-\frac{(a_{n}|x|)^{\frac{1}{c}}}{2}}^{1-\frac{a_{n}|x|}{2}}
(1-s)^{\alpha} s^{\beta} \left( 1+\frac{a_{n}|x|s\left( 1-\frac{3}{2}s \right)+\frac{3}{2}a_{n}^{2}|x|^{2}s^{2}}{1-s} \right)^{\alpha}ds  \nonumber \\
&<&
\int_{1-\frac{(a_{n}|x|)^{\frac{1}{c}}}{2}}^{1-\frac{a_{n}|x|}{2}}
(1-s)^{\alpha} s^{\beta} ds  \nonumber\\
&<&
\frac{(a_{n}|x|)^{\frac{\alpha}{c}}}{2^{\alpha}}\int_{1-\frac{(a_{n}|x|)^{\frac{1}{c}}}{2}}^{1-\frac{a_{n}|x|}{2}} s^{\beta} ds \nonumber \\
&=&
\frac{(a_{n}|x|)^{\frac{\alpha+1}{c}}}{2^{\alpha+1}}(1+o(1))=o(a_{n}),
\end{eqnarray}
\begin{eqnarray}\label{eq3.9}
\int_{1-\frac{(a_{n}|x|)^{\frac{1}{c}}}{2}}^{1} (1-s)^{\alpha-1} s^{\beta+1} ds
< \int_{1-\frac{(a_{n}|x|)^{\frac{1}{c}}}{2}}^{1} (1-s)^{\alpha-1}ds
=\frac{(a_{n}|x|)^{\frac{\alpha}{c}}}{\alpha 2^{\alpha}}=O(a_{n}^{\alpha/\lx{c}}),
\end{eqnarray}
\begin{eqnarray}\label{eq3.10}
\int_{1-\frac{(a_{n}|x|)^{\frac{1}{c}}}{2}}^{1} (1-s)^{\alpha} s^{\beta} ds
\leq \frac{(a_{n}|x|)^{\frac{\alpha}{c}}}{2^{\alpha}(\beta+1)}\left( 1-\left( 1-\frac{(a_{n}|x|)^{\frac{1}{c}}}{2} \right)^{\beta+1} \right)
=o(a_{n})
\end{eqnarray}
and
\begin{eqnarray}\label{eq3.11}
& & a_{n}^{2}\int_{\frac{2}{3}}^{1-\frac{(a_{n}|x|)^{\frac{1}{c}}}{2}} (1-s)^{\alpha-2}s^{\beta+2}\left( 1-\frac{3}{2}s+\frac{3}{2}a_{n}|x|s \right)^{2}ds \nonumber\\
&<& 9a_{n}^{2} \int_{\frac{2}{3}}^{1-\frac{(a_{n}|x|)^{\frac{1}{c}}}{2}} (1-s)^{\alpha-2}ds \nonumber\\
&=& \frac{9a_{n}^{2}}{1-\alpha} \left( \frac{(a_{n}|x|)^{\frac{\alpha-1}{c}}}{2^{\alpha-1}} - 3^{1-\alpha}  \right) \nonumber \\
&<& \frac{9a_{n}^{2+\frac{\alpha-1}{c}}|x|^{\frac{\alpha-1}{c}}}{(1-\alpha)2^{\alpha-1}}
=o(a_{n})
\end{eqnarray}
for large $n$, if $0<\alpha<1$. From \eqref{eq3.4} and \eqref{eq3.9}-\eqref{eq3.11}, it follows that
\begin{eqnarray}\label{eq3.12}
& & \int_{\frac{2}{3}}^{1-\frac{(a_{n}|x|)^{\frac{1}{c}}}{2}}
\left( 1-t_{n}(x)s-\frac{3}{2}a_{n}|x|t_{n}(x)s^{2}  \right)^{\alpha} s^{\beta}ds \nonumber \\
&<& \int_{\frac{2}{3}}^{1-\frac{(a_{n}|x|)^{\frac{1}{c}}}{2}} (1-s)^{\alpha}s^{\beta}ds
+\alpha \int_{\frac{2}{3}}^{1-\frac{(a_{n}|x|)^{\frac{1}{c}}}{2}}(1-s)^{\alpha-1}s^{\beta+1}a_{n}|x|
\left( 1-\frac{3}{2}s+\frac{3}{2}a_{n}|x|s \right)ds \nonumber\\
&=& \int_{\frac{2}{3}}^{1} (1-s)^{\alpha}s^{\beta}ds +\alpha a_{n} |x|\int_{\frac{2}{3}}^{1} (1-s)^{\alpha-1}s^{\beta+1}\left( 1-\frac{3}{2}s \right)ds + o(a_{n})
\end{eqnarray}
and
\begin{eqnarray}\label{eq3.13}
& & \int_{\frac{2}{3}}^{1-\frac{(a_{n}|x|)^{\frac{1}{c}}}{2}}
\left( 1-t_{n}(x)s-\frac{3}{2}a_{n}|x|t_{n}(x)s^{2}  \right)^{\alpha} s^{\beta}ds \nonumber \\
&>& \int_{\frac{2}{3}}^{1-\frac{(a_{n}|x|)^{\frac{1}{c}}}{2}} (1-s)^{\alpha}s^{\beta}ds
+\alpha \int_{\frac{2}{3}}^{1-\frac{(a_{n}|x|)^{\frac{1}{c}}}{2}}(1-s)^{\alpha-1}s^{\beta+1}a_{n}|x|
\left( 1-\frac{3}{2}s+\frac{3}{2}a_{n}|x|s \right)ds \nonumber\\
& & + 2^{1-\alpha}\alpha(\alpha-1)a_{n}^{2}|x|^2\int_{\frac{2}{3}}^{1-\frac{(a_{n}|x|)^{\frac{1}{c}}}{2}}
(1-s)^{\alpha-2}s^{\beta+2}\left(  1-\frac{3}{2}s+\frac{3}{2}a_{n}|x|s \right)^{2}ds\nonumber\\
&=&\int_{\frac{2}{3}}^{1} (1-s)^{\alpha}s^{\beta}ds +\alpha a_{n} |x|\int_{\frac{2}{3}}^{1} (1-s)^{\alpha-1}s^{\beta+1}\left( 1-\frac{3}{2}s \right)ds + o(a_{n})
\end{eqnarray}
for large $n$. Combining with \eqref{eq3.8}, we can get
\begin{eqnarray}\label{eq3.14}
&&\int_{\frac{2}{3}}^{1-\frac{a_{n}|x|}{2}}
\left( 1-t_{n}(x)s-\frac{3}{2}a_{n}|x|t_{n}(x)s^{2}  \right)^{\alpha} s^{\beta}ds\nonumber\\
&=& \int_{\frac{2}{3}}^{1} (1-s)^{\alpha}s^{\beta}ds +\alpha a_{n} |x|\int_{\frac{2}{3}}^{1} (1-s)^{\alpha-1}s^{\beta+1}\left( 1-\frac{3}{2}s \right)ds + o(a_{n}),
\end{eqnarray}
if $0<\alpha<1$.

Now we consider the cases of $\alpha\geq 1$. Since for large $n$ and $\beta>-1$
\begin{eqnarray*}
\int_{1-\frac{a_{n}|x|}{2}}^{1} (1-s)^{\alpha-2}s^{\beta+2}ds
\leq \int_{1-\frac{a_{n}|x|}{2}}^{1} (1-s)^{\alpha-2} ds
=\frac{(a_{n}|x|)^{\alpha-1}}{(\alpha-1)2^{\alpha-1}}
=O(a_{n}^{\alpha-1})
\end{eqnarray*}
and
\begin{eqnarray*}
\int_{1-\frac{a_{n}|x|}{2}}^{1} (1-s)^{\alpha-1}s^{\beta+1}ds
\leq \int_{1-\frac{a_{n}|x|}{2}}^{1} (1-s)^{\alpha-1} ds
=\frac{(a_{n}|x|)^{\alpha}}{\alpha2^{\alpha}}
=o(a_{n})
\end{eqnarray*}
hold with $\alpha>1$, and
\begin{eqnarray*}
\int_{1-\frac{a_{n}|x|}{2}}^{1} (1-s)^{\alpha}s^{\beta}ds
\leq \left(\frac{a_{n}|x|}{2}\right)^{\lx{\alpha}}\int_{1-\frac{a_{n}|x|}{2}}^{1} s^{\beta} ds
=\frac{(a_{n}|x|)^{\alpha+1}}{2^{\alpha+1}}(1+o(1))
=o(a_{n})
\end{eqnarray*}
holds with $\alpha>0$, \eqref{eq3.14} also holds for the cases $\alpha \geq 1$ by using arguments
similar to that of \eqref{eq3.12} and \eqref{eq3.13}.

Similarly, by using \eqref{eq3.5}-\eqref{eq3.7} we have
\begin{eqnarray*}
&&\int_{0}^{\frac{2}{3}}
\left( 1-t_{n}(x)s-\frac{3}{2}a_{n}|x|t_{n}(x)s^{2}  \right)^{\alpha} s^{\beta}ds\\
&=&\int^{\frac{2}{3}}_{0} (1-s)^{\alpha}s^{\beta}ds +\alpha a_{n} |x|\int^{\frac{2}{3}}_{0} (1-s)^{\alpha-1}s^{\beta+1}\left( 1-\frac{3}{2}s \right)ds + o(a_{n})
\end{eqnarray*}
for large $n$, which implies that
\begin{eqnarray}\label{eq3.15}
&&\int_{0}^{1-\frac{a_{n}|x|}{2}}
\left( 1-t_{n}(x)s-\frac{3}{2}a_{n}|x|t_{n}(x)s^{2}  \right)^{\alpha} s^{\beta}ds\nonumber\\
&=& \int_{0}^{1} (1-s)^{\alpha}s^{\beta}ds +\alpha a_{n} |x|\int_{0}^{1} (1-s)^{\alpha-1}s^{\beta+1}\left( 1-\frac{3}{2}s \right)ds + o(a_{n})
\end{eqnarray}
by combining with \eqref{eq3.14}.

 From Lemma \ref{lem1}, \eqref{eq3.15} and the following inequalities
\begin{eqnarray*}
1+\frac{1}{2}s +\frac{3}{8}s^{2}<(1-s)^{-\frac{1}{2}}<1+\frac{1}{2}s +\frac{3}{8}s^{2}+\frac{5\sqrt{2}}{2}s^{3}, 0<s<\frac{1}{2}
\end{eqnarray*}
and
\begin{eqnarray*}
1-(b-1)s<(1-s)^{b-1}<1-(b-1)s+2^{2-b}(b-1)(b-2)s^{2}, 0<s<\frac{1}{2}, b<1,
\end{eqnarray*}
we can get{\small
\begin{eqnarray}\label{eq3.16}
& & \E\left( 1-F(t_{n}(x)(1-Q)^{-\frac{1}{2}}) \right) \nonumber\\
&=&
\frac{\Gamma(a+b)}{\Gamma(a)\Gamma(b)}\int_{0}^{2a_{n}|x|-a_{n}^{2}x^2}
\left( 1-F\left( t_{n}(x)(1-s)^{-\frac{1}{2}} \right)  \right) s^{a-1}(1-s)^{b-1}ds\nonumber\\
&\leq&
\frac{\Gamma(a+b)}{\Gamma(a)\Gamma(b)}\int_{0}^{2a_{n}|x|-a_{n}^{2}x^2}
\left( 1-F\left( t_{n}(x)\left(1+\frac{1}{2}s+\frac{3}{8}s^2  \right) \right)  \right)\nonumber\\
&& \times s^{a-1}\left(1-(b-1)s+2^{2-b}(b-1)(b-2)s^2\right)ds \nonumber\\
&=&
\frac{\Gamma(a+b)}{\Gamma(a)\Gamma(b)}(2a_{n}|x|)^a\int_{0}^{1-\frac{a_{n}|x|}{2}}
\left( 1-F\left( t_{n}(x)\left( 1+a_{n}|x|s+\frac{3}{2}a_{n}^2x^2s^2  \right)  \right) \right)\nonumber\\
& & \times s^{a-1}
\left( 1-2a_{n}|x|(b-1)s+2^{4-b}(b-1)(b-2)a_{n}^2x^2s^2  \right)ds \nonumber\\
&\leq&
\frac{\Gamma(a+b)}{\Gamma(a)\Gamma(b)}(2a_{n}|x|)^a \left( 1-F(t_{n}(x)) \right) \nonumber\\
& & \times
\left(  \int_{0}^{1-\frac{a_{n}|x|}{2}} \left( 1-t_{n}(x)s-\frac{3}{2}a_{n}|x|t_{n}(x)s^2 \right)^{\alpha}
s^{a-1}\left( 1-2a_{n}|x|(b-1)s+2^{4-b}(b-1)(b-2)a_{n}^2x^2s^2  \right)ds \right.\nonumber\\
& &
\left. +A(a_{n}^{-1}|x|^{-1})\int_{0}^{1-\frac{a_{n}|x|}{2}} \left(\left( 1-t_{n}(x)s-\frac{3}{2}a_{n}|x|t_{n}(x)s^2 \right)^{\alpha}
\frac{ \left( 1-t_{n}(x)s-\frac{3}{2}a_{n}|x|t_{n}(x)s^2 \right)^{-\tau}-1}{\tau}\right. \right.\nonumber\\
& &\quad
\left. \left.
+\varepsilon \left( 1-t_{n}(x)s-\frac{3}{2}a_{n}|x|t_{n}(x)s^2 \right)^{\alpha-\tau-\varepsilon}\right)
\times
s^{a-1}\left( 1-2a_{n}|x|(b-1)s+2^{4-b}(b-1)(b-2)a_{n}^2x^2s^2  \right)ds  \right)\nonumber\\
&=&
\frac{\Gamma(a+b)}{\Gamma(a)\Gamma(b)}(2a_{n}|x|)^a \left( 1-F(t_{n}(x)) \right)
\left( \int_{0}^{1}(1-s)^{\alpha}s^{a-1}ds
+A(a_{n}^{-1})|x|^{-\tau}\int_{0}^{1}(1-s)^{\alpha}\frac{(1-s)^{-\tau}-1}{\tau}s^{a-1}ds
\right. \nonumber\\
& & \left.
+a_{n}|x|\int_{0}^{1}(1-s)^{\alpha-1}s^{a}\left( \alpha -2(b-1)-\frac{3}{2}\alpha s + 2(b-1)s \right)ds
+o\left( a_{n}+A(a_{n}^{-1}) \right)
 \right)
\end{eqnarray}}
for large $n$. Similarly,{\small
\begin{eqnarray}\label{eq3.17}
& & \E\left(  1-F\left(  t_{n}(x)(1-Q)^{-\frac{1}{2}} \right) \right) \nonumber\\
&\geq& \frac{\Gamma(a+b)}{\Gamma(a)\Gamma(b)}\int_{0}^{2a_{n}|x|-a_{n}^2x^2}
\left( 1-F\left(t_{n}(x)\left( 1+\frac{1}{2}s+\frac{3}{8}s^2+\frac{5\sqrt{2}}{2}s^3 \right) \right) \right)
s^{a-1}(1-(b-1)s)ds \nonumber\\
&\geq&
\frac{\Gamma(a+b)}{\Gamma(a)\Gamma(b)}(2a_{n}|x|)^a \left( 1-F(t_{n}(x)) \right) \nonumber\\
& & \times
\left(  \int_{0}^{1-\frac{a_{n}|x|}{2}} \left( 1-t_{n}(x)s-\frac{3}{2}a_{n}|x|t_{n}(x)s^2
-20\sqrt{2}a_{n}^2x^2t_{n}(x)s^3 \right)^{\alpha}
s^{a-1}\left( 1-2a_{n}|x|(b-1)s\right)ds \right.\nonumber\\
& &
\left. +A(a_{n}^{-1}|x|^{-1})\int_{0}^{1-\frac{a_{n}|x|}{2}} \left( \left( 1-t_{n}(x)s-\frac{3}{2}a_{n}|x|t_{n}(x)s^2 -20\sqrt{2}a_{n}^2x^2t_{n}(x)s^3 \right)^{\alpha} \right. \right.\nonumber\\
& &
\quad \left. \left.
\times
\frac{ \left( 1-t_{n}(x)s-\frac{3}{2}a_{n}|x|t_{n}(x)s^2
-20\sqrt{2}a_{n}^2x^2t_{n}(x)s^3 \right)^{-\tau}-1}{\tau} \right. \right.\nonumber\\
& &\quad
\left. \left.
-\varepsilon \left( 1-t_{n}(x)s-\frac{3}{2}a_{n}|x|t_{n}(x)s^2 -20\sqrt{2}a_{n}^2x^2t_{n}(x)s^3 \right)^{\alpha-\tau-\varepsilon}\right)
\times
s^{a-1}\left( 1-2a_{n}|x|(b-1)s \right)ds  \right)\nonumber\\
&=&
\frac{\Gamma(a+b)}{\Gamma(a)\Gamma(b)}(2a_{n}|x|)^a \left( 1-F(t_{n}(x)) \right)
\left( \int_{0}^{1}(1-s)^{\alpha}s^{a-1}ds
+A(a_{n}^{-1})|x|^{-\tau}\int_{0}^{1}(1-s)^{\alpha}\frac{(1-s)^{-\tau}-1}{\tau}s^{a-1}ds
\right. \nonumber\\
& & \left.
+a_{n}|x|\int_{0}^{1}(1-s)^{\alpha-1}s^{a}\left( \alpha -2(b-1)-\frac{3}{2}\alpha s + 2(b-1)s \right)ds
+o\left( a_{n}+A(a_{n}^{-1}) \right)
 \right).
\end{eqnarray}}
Combining \eqref{eq3.16} and \eqref{eq3.17}, we can derive \eqref{eq3.1}, which complete the proof. \qed

\begin{lemma}\label{lem3}
Let $(S_{1}, S_{2})$ be a bivariate spherical random vector with almost surely positive random
radius $R$ with df $F$. Assume that $F$ has upper
endpoint $1$ and further $1-F(1-t^{-1})\in 2RV_{-\alpha,\tau}$ with $t>0, \alpha>0, \tau\leq 0$
and auxiliary function $A(t)$.
Then $1-G(1-t^{-1})\in 2RV_{-\alpha-\frac{1}{2},\max(-1,\tau)}$, where $G$ is the df of $S_{1}$.
\end{lemma}

\noindent
{\bf Proof.}
By Corollary 12.1.1 of Berman (1992) and arguments similar to Lemma \ref{lem2}, we can get
\begin{eqnarray}\label{eq3.18}
& & 1-G(1-(t|x|)^{-1})\nonumber\\
&=& \frac{1}{2}\E\left( 1-F\left( (1-(t|x|)^{-1})(1-Q)^{-\frac{1}{2}} \right) \right) \nonumber\\
&=& \frac{(2t^{-1}|x|^{-1})^{\frac{1}{2}}}{2\pi}\left( 1-F\left( 1-(t|x|)^{-1} \right) \right)\nonumber\\
&&\times\left\{  \int_{0}^{1} (1-s)^{\alpha}s^{-\frac{1}{2}}ds +  (t|x|)^{-1}\int_{0}^{1}(1-s)^{\alpha-1}s^{\frac{1}{2}}
\left( \alpha+1-\frac{3}{2}\alpha s -s \right)ds \right. \nonumber\\
& & \left. + A(t)|x|^{\tau}\int_{0}^{1} (1-s)^{\alpha}\frac{(1-s)^{-\tau}-1}{\tau}s^{-\frac{1}{2}}ds
+ o(t^{-1}+A(t))\right\}
\end{eqnarray}
for large $t$ and $x<0$. Here $Q$ is a Beta distributed random variable with parameters ${1}/{2}, {1}/{2}$.
Hence,{\small
\begin{eqnarray}\label{eq3.19}
\frac{1-G(1-(t|x|)^{-1})}{1-G(1-t^{-1})}
&=&|x|^{-\frac{1}{2}-\alpha}\left( 1+\frac{t^{-1}(|x|^{-1}-1)}{c_{\alpha}}\int_{0}^{1}(1-s)^{\alpha-1}s^{\frac{1}{2}}\left( \alpha+1-\frac{3}{2}\alpha s -s \right)ds \right. \nonumber\\
& & \left. + A(t)\left(  \frac{|x|^{\tau}-1}{\tau} + \frac{|x|^{\tau}-1}{c_{\alpha}}
\int_{0}^{1}(1-s)^{\alpha}\frac{(1-s)^{-\tau}-1}{\tau}s^{-\frac{1}{2}}ds \right) +o(t^{-1}+A(t)) \right) \nonumber\\
& &
\end{eqnarray}}
for large $t$, where $c_{\alpha}=\int_{0}^{1}(1-s)^{\alpha}s^{-\frac{1}{2}}ds$. So $1-G(1-t^{-1})\in 2RV_{-\alpha-\frac{1}{2},\max(-1,\tau)}$.
\qed

\section{Proofs of main results}
\label{sec4}

The aim of this section is to prove our main results.

\noindent
{\bf Proof of Theorem \ref{th1}.}
Note that $\arccos z =\sqrt{2(1-z)}\left( 1+\frac{1-z}{12}+o(1-z) \right)$ for $z\to 1$.
By \eqref{eq1.4} and the Taylor expansion with the Lagrange remainder, we have{\small
\begin{eqnarray}\label{eq4.1}
\frac{\frac{t_{n}(y)}{t_{n}(x)}-\rho_{n}}{\sqrt{a_{n}}\sin \tau_{n}}
&=& \frac{1-\rho_{n}+\frac{a_{n}(y-x)}{1+a_{n}x}}{\sqrt{a_{n}}\left( \sqrt{2(1-\rho_{n})}
\left( 1+\frac{1}{12}(1-\rho_{n})+o(1-\rho_{n}) \right) -\frac{(\sqrt{2(1-\rho_{n})})^3}{6}(1+o(1)) \right)} \nonumber\\
&=& \lambda+ \frac{y-x}{2\lambda} +\gamma\left( 1-\frac{y-x}{2\lambda^2} \right)c_{n}
+\left( \frac{\lambda^3}{2} +\frac{(y-x)\lambda}{4} -\frac{x(y-x)}{2\lambda} \right) a_{n} +o(a_{n}+c_{n})
\end{eqnarray}}
for large $n$, where $\tau_{n}=\arccos \rho_{n}$. Then for
 $\beta_{n}(x,y)=\arctan\left( \frac{\frac{y}{x}-\rho_{n}}{\sin \tau_{n}} \right)$, we can get
\begin{eqnarray}\label{eq4.2}
& & \frac{\sin \beta_{n}\left( t_{n}(x), t_{n}(y) \right)}{\sqrt{a_{n}}} \nonumber \\
&=& \frac{\beta_{n}\left( t_{n}(x), t_{n}(y) \right)}{\sqrt{a_{n}}}
-\frac{\beta_{n}^{3}\left( t_{n}(x) ,t_{n}(y) \right)}{6\sqrt{a_{n}}}(1+o(1)) \nonumber \\
&=& \lambda +\frac{y-x}{2\lambda} +\gamma\left( 1-\frac{y-x}{2\lambda^2} \right)c_{n}\nonumber \\
&&-\frac{a_{n}}{2}\left( \lambda(y-x)+\frac{x(y-x)}{\lambda} + \frac{3(y-x)^2}{4\lambda}
+ \frac{(y-x)^3}{8\lambda^3}  \right) +o(a_{n}+c_{n})
\end{eqnarray}
for large $n$.

First we consider the case of $\lambda+\frac{y-x}{2\lambda}\geq 0$ with $\lambda^2+x+y+\frac{(x-y)^2}{4\lambda^2}<0$. Combining \eqref{eq4.2} and Taylor's Theorem,
we can get
\begin{eqnarray*}
& & \int_{\frac{\sin^2\beta_{n}(t_{n}(x),t_{n}(y))}{2a_{n}|x|}}^{1-\frac{a_{n}|x|}{2}}
\left( 1-t_{n}(x)s-\frac{3}{2}a_{n}|x|t_{n}(x)s^2 \right)^{\alpha}s^{\beta}ds \nonumber\\
&=& \int_{\frac{\left( \lambda+\frac{y-x}{2\lambda} \right)^2}{2|x|}}^{1-\frac{a_{n}|x|}{2}}
\left( 1-t_{n}(x)s-\frac{3}{2}a_{n}|x|t_{n}(x)s^2 \right)^{\alpha}s^{\beta}ds -
\left( \frac{\sin^2\beta_{n}(t_{n}(x),t_{n}(y))}{2a_{n}|x|}-
\frac{\left( \lambda+\frac{y-x}{2\lambda} \right)^2}{2|x|}\right) \nonumber\\
& & \times \left( 1-\frac{1+a_{n}x}{2|x|}\left( \lambda+\frac{y-x}{2\lambda} \right)^2
-\frac{3}{2}a_{n}|x|\frac{\left( \lambda+\frac{y-x}{2\lambda} \right)^4}{4x^2}(1+a_{n}x) \right)^{\alpha}
\frac{\left( \lambda+\frac{y-x}{2\lambda} \right)^{2\beta}}{(2|x|)^{\beta}}(1+o(1)) \nonumber\\
&=& \int_{\frac{\left( \lambda+\frac{y-x}{2\lambda} \right)^2}{2|x|}}^{1} (1-s)^{\alpha}s^{\beta}ds
+ \alpha a_{n}|x|\int_{\frac{\left( \lambda+\frac{y-x}{2\lambda} \right)^2}{2|x|}}^{1} (1-s)^{\alpha-1}
s^{\beta+1}\left( 1-\frac{3}{2}s \right)ds \nonumber\\
& & -\frac{\left( \lambda+\frac{y-x}{2\lambda} \right)^{2\beta+1}}{2^\beta|x|^{\beta+1}}
\left( 1-\frac{\left( \lambda+\frac{y-x}{2\lambda} \right)^2}{2|x|} \right)^{\alpha}
\left\{ c_{n}\gamma\left( 1-\frac{y-x}{2\lambda^2} \right)\right.\\
&&\left. -\frac{a_{n}}{2}\left( \lambda(y-x)+\frac{x(y-x)}{\lambda}+\frac{3(y-x)^2}{4\lambda} +\frac{(y-x)^3}{8\lambda^3}  \right) \right\} +o(a_{n}+c_{n})
\end{eqnarray*}
for large $n$ and $\beta>-1$, since
\begin{eqnarray*}
& & \int_{\frac{\left( \lambda+\frac{y-x}{2\lambda} \right)^2}{2|x|}}^{1-\frac{a_{n}|x|}{2}}
\left( 1-t_{n}(x)s-\frac{3}{2}a_{n}|x|t_{n}(x)s^2 \right)^{\alpha}s^{\beta}ds \nonumber\\
&=&\int_{\frac{\left( \lambda+\frac{y-x}{2\lambda} \right)^2}{2|x|}}^{1} (1-s)^{\alpha}s^{\beta}ds
+\alpha a_{n}|x|\int_{\frac{\left( \lambda+\frac{y-x}{2\lambda} \right)^2}{2|x|}}^{1} (1-s)^{\alpha-1} s^{\beta+1}
\left( 1-\frac{3}{2}s \right)ds +o(a_n+c_n)
\end{eqnarray*}
holds for large $n$ by arguments similar with \eqref{eq3.15}.

With the same arguments as Lemma \ref{lem2}, we can get{\small
\begin{eqnarray*}
& & \frac{1}{2\pi}\int_{\beta_{n}\left(t_{n}(x),t_{n}(y)\right)}^{\frac{\pi}{2}}\left( 1-F\left( t_{n}(x)/\cos \alpha \right) \right)d\alpha \nonumber\\
&=& \frac{1}{4\pi} \int_{\sin^2 \beta_{n}\left(t_{n}(x), t_{n}(y)\right)}^{2a_{n}|x|-a_{n}^2x^2}
\left( 1-F\left( t_{n}(x)(1-s)^{-\frac{1}{2}}  \right)  \right) s^{-\frac{1}{2}} (1-s)^{-\frac{1}{2}}ds \nonumber\\
&=& \frac{(2a_{n}|x|)^{\frac{1}{2}}}{4\pi}\left( 1-F(t_{n}(x)) \right)
\left\{   \int_{\frac{\left( \lambda+\frac{y-x}{2\lambda} \right)^2}{2|x|}}^{1} (1-s)^{\alpha}s^{-\frac{1}{2}} ds
+A(a_{n}^{-1})|x|^{-\tau}\int_{\frac{\left( \lambda+\frac{y-x}{2\lambda} \right)^2}{2|x|}}^{1} (1-s)^{\alpha}\frac{(1-s)^{-\tau}-1}{\tau}s^{-\frac{1}{2}}ds \right. \nonumber\\
& &
\left.  +a_{n}|x|\int^{1}_{\frac{\left( \lambda+\frac{y-x}{2\lambda} \right)^{2}}{2|x|}} (1-s)^{\alpha-1}s^{\frac{1}{2}}\left(  \alpha+1-s-\frac{3}{2}\alpha s \right) ds
+ o(a_{n}+c_{n}+A(a_{n}^{-1})) \right. \nonumber\\
& &
\left. -\left( \frac{2}{|x|} \right)^{\frac{1}{2}} \left( 1-\frac{\left( \lambda+\frac{y-x}{2\lambda} \right)^2}{2|x|} \right)^{\alpha}
\left( c_{n}\gamma\left( 1-\frac{y-x}{2\lambda^2} \right) -\frac{a_{n}}{2}\left( \lambda(y-x)+\frac{x(y-x)}{\lambda}+\frac{3(y-x)^2}{4\lambda} +\frac{(y-x)^3}{8\lambda^3}  \right) \right)
\right\}
\end{eqnarray*}}
for large $n$, which implies that
\begin{eqnarray}\label{eq4.3}
& & \frac{n}{2\pi}\int_{\beta_{n}(t_{n}(x),t_{n}(y))}^{\frac{\pi}{2}}\left( 1-F(t_{n}(x)/\cos \alpha) \right)d\alpha\nonumber\\
&=& |x|^{\alpha+\frac{1}{2}} \left\{ \frac{1}{c_{\alpha}}\int^{1}_{\frac{\lambda+\frac{y-x}{2\lambda}}{\sqrt{2|x|}}}
(1-s^2)^{\alpha}ds
+\frac{A(a_{n}^{-1})}{c_{\alpha}}\frac{|x|^{-\tau}-1}{\tau} \int^{1}_{\frac{\lambda+\frac{y-x}{2\lambda}}{\sqrt{2|x|}}}
(1-s^2)^{\alpha}ds \right. \nonumber\\
& & \left. -\frac{a_{n}}{c_{\alpha}^2}\int_{-1}^{1}(1-s^2)^{\alpha-1}s^2\left( \alpha+1-s^2-\frac{3}{2}\alpha s^2 \right)ds \int^{1}_{\frac{\lambda+\frac{y-x}{2\lambda}}{\sqrt{2|x|}}}
(1-s^2)^{\alpha}ds \right. \nonumber\\
& & -\frac{A(a_{n}^{-1})}{c_{\alpha}^2}\int_{-1}^{1}(1-s^2)^{\alpha}\frac{(1-s^2)^{-\tau}-1}{\tau}ds
\int^{1}_{\frac{\lambda+\frac{y-x}{2\lambda}}{\sqrt{2|x|}}}
(1-s^2)^{\alpha}ds\nonumber\\
&& \left.+\frac{A(a_{n}^{-1})|x|^{-\tau}}{c_{\alpha}}\int_{\frac{\lambda+\frac{y-x}{2\lambda}}{\sqrt{2|x|}}}^{1}
(1-s^2)^{\alpha}\frac{(1-s^2)^{-\tau}-1}{\tau}ds \right.  \nonumber\\
& &  \left. + \frac{a_{n}|x|}{c_{\alpha}}\int_{\frac{\lambda+\frac{y-x}{2\lambda}}{\sqrt{2|x|}}}^{1}
(1-s^2)^{\alpha-1}s^2\left( \alpha+1-s^2-\frac{3}{2}\alpha s^2  \right)ds
+ o(a_{n}+c_{n}+A(a_{n}^{-1})) \right. \nonumber\\
& & -\frac{1}{c_{\alpha}\sqrt{2|x|}}\left( 1-\frac{\left( \lambda+\frac{y-x}{2\lambda} \right)^2}{2|x|} \right)^{\alpha}
\left( c_{n}\gamma\left( 1-\frac{y-x}{2\lambda^2} \right)\right.\nonumber\\
&&\left.\left.-\frac{a_{n}}{2}\left( \lambda(y-x)+\frac{x(y-x)}{\lambda}+\frac{3(y-x)^2}{4\lambda} +\frac{(y-x)^3}{8\lambda^3}  \right) \right)
\right\} \nonumber\\
&:=& M_{\alpha+\frac{1}{2},\lambda}(x,y).
\end{eqnarray}

Similarly, for the case of $\lambda+\frac{y-x}{2\lambda}< 0$ with $\lambda^2+x+y+\frac{(x-y)^2}{4\lambda^2}<0$, we have
\begin{eqnarray*}
& & \frac{n}{2\pi}\int_{\beta_{n}(t_{n}(x),t_{n}(y))}^{\frac{\pi}{2}}\left( 1-F(t_{n}(x)/\cos \alpha) \right)d\alpha \nonumber\\
&=& \frac{n}{4\pi}\int_{0}^{1} \left( 1-F\left( t_{n}(x)(1-s)^{-\frac{1}{2}} \right)  \right)s^{-\frac{1}{2}}
(1-s)^{-\frac{1}{2}}ds \\
& &
+\frac{n}{4\pi}\int_{0}^{\sin^2\beta_{n}(t_{n}(x),t_{n}(y))} \left( 1-F\left( t_{n}(x)(1-s)^{-\frac{1}{2}} \right)  \right)s^{-\frac{1}{2}}
(1-s)^{-\frac{1}{2}}ds \\
&=& M_{\alpha+\frac{1}{2},\lambda}(x,y)
\end{eqnarray*}
for large $n$, which is consistent with the case of $\lambda+\frac{y-x}{2\lambda}\geq 0$.

Define $\tilde{\beta}_{n}(t_{n}(x),t_{n}(y))=\tau_{n}-\beta_{n}(t_{n}(x),t_{n}(y))$. From the second-order condition \eqref{eq1.4} and \eqref{eq4.2}, it follows that
\begin{eqnarray}\label{eq4.4}
& & \frac{\tilde{\beta}_{n}(t_{n}(x),t_{n}(y))}{\sqrt{a_{n}}} \nonumber\\
&=& \lambda+\frac{x-y}{2\lambda}+ \gamma\left(  1-\frac{x-y}{2\lambda^2} \right)c_{n}\nonumber\\
&&+\frac{a_{n}}{2}\left( \frac{\lambda^2}{3} +\frac{\lambda(y-x)}{2}+\frac{x(y-x)}{\lambda}+\frac{(y-x)^2}{2\lambda}
+\frac{(y-x)^3}{12\lambda^3}  \right) +o(a_{n}+c_{n})
\end{eqnarray}
for large $n$, which implies that
\begin{eqnarray}\label{eq4.5}
\frac{\sin \tilde{\beta}_{n}(t_{n}(x),t_{n}(y))}{\sqrt{a_{n}}}
&=& \lambda+\frac{x-y}{2\lambda} +\gamma\left(  1-\frac{x-y}{2\lambda^2} \right)c_{n}\nonumber\\
&&+ \frac{a_{n}}{2}\left( \lambda(y-x) +\frac{x(y-x)}{\lambda} +\frac{(y-x)^2}{4\lambda}+\frac{(y-x)^3}{8\lambda^3} \right)
\end{eqnarray}
for large $n$. Using the similar arguments as above, we can get
\begin{eqnarray}\label{eq4.6}
& & \frac{n}{2\pi}\int_{\tilde{\beta}_{n}(t_{n}(x),t_{n}(y))}^{\frac{\pi}{2}}\left( 1-F(t_{n}(y)/\cos \alpha) \right)d\alpha\nonumber\\
&=& |y|^{\alpha+\frac{1}{2}} \left\{ \frac{1}{c_{\alpha}}\int^{1}_{\frac{\lambda+\frac{x-y}{2\lambda}}{\sqrt{2|y|}}}
(1-s^2)^{\alpha}ds
+\frac{A(a_{n}^{-1})}{c_{\alpha}}\frac{|y|^{-\tau}-1}{\tau} \int^{1}_{\frac{\lambda+\frac{x-y}{2\lambda}}{\sqrt{2|y|}}}
(1-s^2)^{\alpha}ds \right. \nonumber\\
& & \left. -\frac{a_{n}}{c_{\alpha}^2}\int_{-1}^{1}(1-s^2)^{\alpha-1}s^2\left( \alpha+1-s^2-\frac{3}{2}\alpha s^2 \right)ds \int^{1}_{\frac{\lambda+\frac{x-y}{2\lambda}}{\sqrt{2|y|}}}
(1-s^2)^{\alpha}ds \right. \nonumber\\
& & -\frac{A(a_{n}^{-1})}{c_{\alpha}^2}\int_{-1}^{1}(1-s^2)^{\alpha}\frac{(1-s^2)^{-\tau}-1}{\tau}ds
\int^{1}_{\frac{\lambda+\frac{x-y}{2\lambda}}{\sqrt{2|y|}}}
(1-s^2)^{\alpha}ds\nonumber\\
&&+\frac{A(a_{n}^{-1})|y|^{-\tau}}{c_{\alpha}}\int_{\frac{\lambda+\frac{x-y}{2\lambda}}{\sqrt{2|y|}}}^{1}
(1-s^2)^{\alpha}\frac{(1-s^2)^{-\tau}-1}{\tau}ds   \nonumber\\
& &   + \frac{a_{n}|y|}{c_{\alpha}}\int_{\frac{\lambda+\frac{x-y}{2\lambda}}{\sqrt{2|y|}}}^{1}
(1-s^2)^{\alpha-1}s^2\left( \alpha+1-s^2-\frac{3}{2}\alpha s^2  \right)ds
+ o(a_{n}+c_{n}+A(a_{n}^{-1})) \nonumber\\
& & -\frac{1}{c_{\alpha}\sqrt{2|y|}}\left( 1-\lx{\frac{\left(\lambda+\frac{x-y}{2\lambda}\right)^2}{2|y|}} \right)^{\alpha}
\left( c_{n}\gamma\left( 1-\frac{x-y}{2\lambda^2} \right)\right. \nonumber\\
&&\left.\left.+\frac{a_{n}}{2}\left( \lambda(y-x)+\frac{x(y-x)}{\lambda}+\lx{\frac{(y-x)^2}{4\lambda}} +\frac{(y-x)^3}{8\lambda^3}  \right) \right)
\right\}
\end{eqnarray}
for large $n$.

Since $\lambda^2 + x+y +\frac{(x-y)^2}{4\lambda^2}<0$ ensures that both $\frac{\left( \lambda+\frac{x-y}{2\lambda} \right)^2}{2|y|}<1$ and $\frac{\left( \lambda+\frac{y-x}{2\lambda} \right)^2}{2|x|}<1$ hold,
and
\begin{eqnarray*}
\P(\xi_{n1}>t_{n}(x), \eta_{n1}>t_{n}(y))
&=&\frac{1}{2\pi}\int^{\frac{\pi}{2}}_{\max(\min(\beta_{n}(t_{n}(x),t_{n}(y)),\frac{\pi}{2}),-\frac{\pi}{2})}
\left( 1-F(t_{n}(x)/\cos \alpha)  \right)d\alpha \\
& & +\frac{1}{2\pi}\int^{\frac{\pi}{2}}_{\max(\min(\tilde{\beta}_{n}(t_{n}(x),t_{n}(y)),\frac{\pi}{2}),-\frac{\pi}{2})}
\left( 1-F(t_{n}(y)/\cos \alpha)  \right)d\alpha,
\end{eqnarray*}
c.f., Hashorva (2006), we have
\begin{eqnarray*}
& & -n\left( 1-\P(\xi_{n1} \leq t_{n}(x), \eta_{n1} \leq t_{n}(y))\right) + |x|^{\alpha+\frac{1}{2}}\psi_{\alpha}
\left( \frac{ \lambda+\frac{y-x}{2\lambda} }{\sqrt{2|x|}} \right)
+ |y|^{\alpha+\frac{1}{2}}\psi_{\alpha}
\left( \frac{ \lambda+\frac{x-y}{2\lambda} }{\sqrt{2|y|}} \right) \\
&=& -n\left( 1-G(t_{n}(x)) \right) +|x|^{\alpha+\frac{1}{2}}
-n\left(  1-G(t_{n}(y)) \right) +|y|^{\alpha+\frac{1}{2}} \\
& & +\frac{n}{2\pi}\int_{\tilde{\beta}_{n}(t_{n}(x),t_{n}(y))}^{\frac{\pi}{2}}\left( 1-F(t_{n}(y)/\cos \alpha) \right)d\alpha -|y|^{\alpha+\frac{1}{2}}\left(  1- \psi_{\alpha}
\left( \frac{ \lambda+\frac{x-y}{2\lambda} }{\sqrt{2|y|}} \right) \right) \\
& & + \frac{n}{2\pi}\int_{\beta_{n}(t_{n}(x),t_{n}(y))}^{\frac{\pi}{2}}\left( 1-F(t_{n}(x)/\cos \alpha) \right)d\alpha -|x|^{\alpha+\frac{1}{2}}\left(  1- \psi_{\alpha}
\left( \frac{ \lambda+\frac{y-x}{2\lambda} }{\sqrt{2|x|}} \right) \right)\\
&=& Q_{\alpha+\frac{1}{2},\lambda}(x,y)
\end{eqnarray*}
for large $n$, where $Q_{\alpha+\frac{1}{2},\lambda}(x,y)$ is the one given by Theorem \ref{th1}. Hence,{\small
\begin{eqnarray*}
& & \P\left(  M_{n1}\leq t_{n}(x), M_{n2}\leq t_{n}(y) \right) -H_{\alpha+\frac{1}{2},\lambda}(x,y)\\
&=& {\P}^{n}\left(  \xi_{n1} \leq t_{n}(x), \eta_{n1} \leq t_{n}(y) \right) -H_{\alpha+\frac{1}{2},\lambda}(x,y)\\
&=&  H_{\alpha+\frac{1}{2},\lambda}(x,y) \left( \exp\left(  n\log \P\left( \xi_{n1} \leq t_{n}(x), \eta_{n1} \leq t_{n}(y) \right)\lx{+}
|x|^{\alpha+\frac{1}{2}}\psi_{\alpha}
\left( \frac{ \lambda+\frac{y-x}{2\lambda} }{\sqrt{2|x|}} \right)
+ |y|^{\alpha+\frac{1}{2}}\psi_{\alpha}
\left( \frac{ \lambda+\frac{x-y}{2\lambda} }{\sqrt{2|y|}} \right) \right)-1 \right)\\
&=& H_{\alpha+\frac{1}{2},\lambda}(x,y) \left( n\log \P(\xi_{n1} \leq t_{n}(x), \eta_{n1} \leq t_{n}(y))
+ |x|^{\alpha+\frac{1}{2}}\psi_{\alpha}
\left( \frac{ \lambda+\frac{y-x}{2\lambda} }{\sqrt{2|x|}} \right)
+ |y|^{\alpha+\frac{1}{2}}\psi_{\alpha}
\left( \frac{ \lambda+\frac{x-y}{2\lambda} }{\sqrt{2|y|}} \right) \right.\\
& & \left. + \frac{1}{2}\left(  n\log \P(\xi_{n1} \leq t_{n}(x), \eta_{n1} \leq t_{n}(y))
+ |x|^{\alpha+\frac{1}{2}}\psi_{\alpha}
\left( \frac{ \lambda+\frac{y-x}{2\lambda} }{\sqrt{2|x|}} \right)
+ |y|^{\alpha+\frac{1}{2}}\psi_{\alpha}
\left( \frac{ \lambda+\frac{x-y}{2\lambda} }{\sqrt{2|y|}} \right)   \right)^2(1+o(1))  \right)\\
&=& H_{\alpha+\frac{1}{2},\lambda}(x,y) \left(   -n\left( 1-\P(\xi_{n1} \leq t_{n}(x), \eta_{n1} \leq t_{n}(y))\right) + |x|^{\alpha+\frac{1}{2}}\psi_{\alpha}
\left( \frac{ \lambda+\frac{y-x}{2\lambda} }{\sqrt{2|x|}} \right)
+ |y|^{\alpha+\frac{1}{2}}\psi_{\alpha}
\left( \frac{ \lambda+\frac{x-y}{2\lambda} }{\sqrt{2|y|}} \right)\right.\\
& & \left.  -\frac{n}{2}\left( 1- \P(\xi_{n1} \leq t_{n}(x), \eta_{n1} \leq t_{n}(y)) \right)^2
-\frac{n}{3}\left( 1- \P(\xi_{n1} \leq t_{n}(x), \eta_{n1} \leq t_{n}(y)) \right)^3(1+o(1)) \right. \\
& &  \left.  + \frac{1}{2}\left(  n\log \P(\xi_{n1} \leq t_{n}(x), \eta_{n1} \leq t_{n}(y))
+ |x|^{\alpha+\frac{1}{2}}\psi_{\alpha}
\left( \frac{ \lambda+\frac{y-x}{2\lambda} }{\sqrt{2|x|}} \right)
+ |y|^{\alpha+\frac{1}{2}}\psi_{\alpha}
\left( \frac{ \lambda+\frac{x-y}{2\lambda} }{\sqrt{2|y|}} \right)   \right)^2(1+o(1))  \right)\\
&=&  H_{\alpha+\frac{1}{2},\lambda}(x,y) \Big(  Q_{\alpha+\frac{1}{2},\lambda}(x,y)
-\frac{1}{2}n^{-1}\left(|x|^{\alpha+\frac{1}{2}}\psi_{\alpha}
\left( \frac{ \lambda+\frac{y-x}{2\lambda} }{\sqrt{2|x|}} \right)
+ |y|^{\alpha+\frac{1}{2}}\psi_{\alpha}
\left( \frac{ \lambda+\frac{x-y}{2\lambda} }{\sqrt{2|y|}} \right)  \right)^{2} \\
& &
+ o\left( a_{n}+c_{n}+A(a_{n}^{-1})+n^{-1} \right)   \Big)
\end{eqnarray*}}
for large $n$. \qed

Because the proofs of Theorem \ref{th2}-\ref{th7} are similar, we only prove Theorem \ref{th5}-\ref{th7} below,
and omit the proofs of Theorem \ref{th2}-\ref{th4}.

\noindent
{\bf Proof of Theorem \ref{th5}.}
For large $n$, we have{\small
\begin{eqnarray}\label{eq4.7}
& & \frac{n}{2\pi}\int_{\min\left( \beta_{n}(t_{n}(x), t_{n}(y)),\frac{\pi}{2} \right)}^{\frac{\pi}{2}}
\left( 1-F(t_{n}(x)/\cos \alpha) \right) d\alpha  \nonumber\\
&=& \frac{n}{4\pi} \int_{\lx{\min\left( \sin^2 \beta_{n}(t_{n}(x),t_{n}(y)) ,1\right)}}^{1}
\left( 1-F(t_{n}(x)(1-s)^{-\frac{1}{2}}) \right)s^{-\frac{1}{2}}(1-s)^{-\frac{1}{2}}ds \nonumber\\
&=& \frac{n}{4\pi}(2a_{n}|x|)^{\frac{1}{2}} \int_{\min \left( \frac{\sin^2
\beta_{n}(t_n(x),t_{n}(y))}{2a_{n}|x|}, 1-\frac{a_{n}|x|}{2}  \right)}^{1-\frac{a_{n}|x|}{2}}
\left( 1-F\left( t_{n}(x)(1-2a_{n}|x|s)^{-\frac{1}{2}} \right)  \right) s^{-\frac{1}{2}}(1-2a_{n}|x|s)^{-\frac{1}{2}}ds
\nonumber \\
&=& 0
\end{eqnarray}}
if $\lambda+\frac{y-x}{2\lambda}>\sqrt{2|x|}$; and if $\lambda+\frac{y-x}{2\lambda}<-\sqrt{2|x|}$, by using
\eqref{eq3.19} we have{\small
\begin{eqnarray}\label{eq4.8}
& & \frac{n}{2\pi} \int_{\max \left( \beta_{n}(t_{n}(x),t_{n}(y)) , -\frac{\pi}{2}\right)}^{\frac{\pi}{2}}
\left( 1-F(t_{n}(x)/\cos \alpha) \right) d\alpha  \nonumber\\
&=&
\frac{n}{4\pi}\left( \int_{0}^{1} \left( 1-F\left( t_{n}(x)(1-s)^{-\frac{1}{2}} \right) \right)s^{-\frac{1}{2}}(1-s)^{-\frac{1}{2}}ds  \right. \nonumber\\
& & \left.
+\int_{0}^{\min\left( \sin^2\beta_{n}(t_{n}(x),t_{n}(y)),1  \right)} \left( 1-F\left( t_{n}(x)(1-s)^{-\frac{1}{2}} \right) \right)s^{-\frac{1}{2}}(1-s)^{-\frac{1}{2}}ds \right)  \nonumber\\
&=& \frac{n}{2\pi} (2a_{n}|x|)^{\frac{1}{2}} \int_{0}^{1-\frac{a_{n}|x|}{2}} \left( 1-F\left( t_{n}(x)(1-2a_{n}|x|s)^{-\frac{1}{2}} \right)  \right)s^{-\frac{1}{2}}(1-2a_{n}|x|s)^{-\frac{1}{2}}ds  \nonumber\\
&=& |x|^{\frac{1}{2}+\alpha} \left\{
1+\frac{a_{n}(|x|-1)}{c_{\alpha}}\int_{0}^{1}(1-s)^{\alpha-1}s^{\frac{1}{2}}\left( \alpha+1-s-\frac{3}{2}\alpha s \right) ds \right. \nonumber\\
& &  \left.
+A(a_{n}^{-1})\left( \frac{|x|^{-\tau}-1}{\tau} +\frac{|x|^{-\tau}-1}{c_{\alpha}} \int_{0}^{1}(1-s)^{\alpha}\frac{(1-s)^{-\tau}-1}{\tau}s^{-\frac{1}{2}}ds  \right)
+o\left( a_{n}+A(a_{n}^{-1}) \right)  \right\}
\end{eqnarray}}
for large $n$.

Similarly, for large $n$ we can get
\begin{eqnarray}\label{eq4.9}
\frac{n}{2\pi}\int_{\min\left( \tilde{\beta}_{n}(t_{n}(x), t_{n}(y)),\frac{\pi}{2} \right)}^{\frac{\pi}{2}}
\left( 1-F(t_{n}(y)/\cos \alpha) \right) d\alpha =0
\end{eqnarray}
if $\lambda+\frac{x-y}{2\lambda}>\sqrt{2|y|}$; and{\small
\begin{eqnarray}\label{eq4.10}
& &\frac{n}{2\pi}\int_{\max \left( \tilde{\beta}_{n}(t_{n}(x), t_{n}(y)), -\frac{\pi}{2} \right)}^{\frac{\pi}{2}}
\left( 1-F(t_{n}(y)/\cos \alpha) \right) d\alpha  \nonumber\\
&=& |y|^{\frac{1}{2}+\alpha} \left\{
1+\frac{a_{n}(|y|-1)}{c_{\alpha}}\int_{0}^{1}(1-s)^{\alpha-1}s^{\frac{1}{2}}\left( \alpha+1-s-\frac{3}{2}\alpha s \right) ds \right. \nonumber\\
& &  \left.
+A(a_{n}^{-1})\left( \frac{|y|^{-\tau}-1}{\tau} +\frac{|y|^{-\tau}-1}{c_{\alpha}} \int_{0}^{1}(1-s)^{\alpha}\frac{(1-s)^{-\tau}-1}{\tau}s^{-\frac{1}{2}}ds  \right)
+o\left( a_{n}+A(a_{n}^{-1}) \right)  \right\}
\end{eqnarray}}
if $\lambda+\frac{x-y}{2\lambda}<-\sqrt{2|y|}$.

Now we only prove case (iii), since the proofs of the rest cases are similar. From \eqref{eq3.19}, \eqref{eq4.7}
and \eqref{eq4.10}, it follows that
\begin{eqnarray*}
& & -n\left( 1-\P(\xi_{n1} \leq t_{n}(x), \eta_{n1} \leq t_{n}(y))\right) + |x|^{\alpha+\frac{1}{2}}\psi_{\alpha}
\left( \frac{ \lambda+\frac{y-x}{2\lambda} }{\sqrt{2|x|}} \right)
+ |y|^{\alpha+\frac{1}{2}}\psi_{\alpha}
\left( \frac{ \lambda+\frac{x-y}{2\lambda} }{\sqrt{2|y|}} \right) \\
&=& -n\left( 1-G(t_{n}(x)) \right) +|x|^{\alpha+\frac{1}{2}}
-n\left(  1-G(t_{n}(y)) \right) +|y|^{\alpha+\frac{1}{2}} \\
& & +\frac{n}{2\pi}\int_{\max\left(\tilde{\beta}_{n}(t_{n}(x),t_{n}(y)), -\frac{\pi}{2} \right)}^{\frac{\pi}{2}}\left( 1-F(t_{n}(y)/\cos \alpha) \right)d\alpha -|y|^{\alpha+\frac{1}{2}} \\
& & + \frac{n}{2\pi}\int_{\min\left(\beta_{n}(t_{n}(x),t_{n}(y)),\frac{\pi}{2}\right)}^{\frac{\pi}{2}}\left( 1-F(t_{n}(x)/\cos \alpha) \right)d\alpha \\
&=&-|x|^{\frac{1}{2}+\alpha} \left(
\frac{a_{n}(|x|-1)}{c_{\alpha}}\int_{0}^{1}(1-s)^{\alpha-1}s^{\frac{1}{2}}\left( \alpha+1-s-\frac{3}{2}\alpha s \right) ds \right. \nonumber\\
& &  \left.
+A(a_{n}^{-1})\left( \frac{|x|^{-\tau}-1}{\tau} +\frac{|x|^{-\tau}-1}{c_{\alpha}} \int_{0}^{1}(1-s)^{\alpha}\frac{(1-s)^{-\tau}-1}{\tau}s^{-\frac{1}{2}}ds  \right)
+o\left( a_{n}+A(a_{n}^{-1}) \right)  \right)
\end{eqnarray*}
for large $n$, then \eqref{eq2.7} holds, which complete the proof. \qed

\noindent
{\bf Proof of Theorem \ref{th6}.}
(i) For the case of $\rho_{n}\equiv 1$. Note that \[\P\left( M_{n1}\leq 1+a_{n}x, M_{n2}\leq 1+a_{n}y \right)
=G^{n}\left( 1+a_{n}\min(x,y) \right)\]and \[H_{\alpha,0}(x,y)=\Psi_{\alpha+\frac{1}{2}}\left( \min(x,y) \right)
=\exp\left( -|\min(x,y)|^{\alpha+\frac{1}{2}}  \right).\]  It follows from \eqref{eq3.19} that \eqref{eq2.9} holds.

(ii) If $\rho_{n}\in (0,1)$ such that $\frac{a_{n}^{{5}/{3}}}{1-\rho_{n}}\to 0$ as $n\to \infty$, which
implies that
\begin{eqnarray*}
\frac{\frac{t_{n}(y)}{t_{n}(x)}-\rho_{n}}{\sqrt{a_{n}} \sin \tau_{n}}
= \lambda_{n}+\frac{y-x}{2\lambda_{n}} -\frac{(y-x)x}{2\lambda_{n}}a_{n} +o\left( \frac{a_{n}}{\lambda_{n}} \right).
\end{eqnarray*}
Hence,
\begin{eqnarray*}
\frac{\sin \beta_{n}(t_{n}(x),t_{n}(y))}{\sqrt{a_{n}}}
= \lambda_{n} +\frac{y-x}{2\lambda_{n}} - \frac{(y-x)^3a_{n}}{16\lambda_{n}^{3}} + o\left( \frac{a_{n}}{\lambda_n^3} \right).
\end{eqnarray*}
Similarly,
\begin{eqnarray*}
\frac{\sin \tilde{\beta}_{n}(t_{n}(x),t_{n}(y))}{\sqrt{a_{n}}}
= \lambda_{n} +\frac{x-y}{2\lambda_{n}} + \frac{(y-x)^3a_{n}}{16\lambda_{n}^{3}} + o\left( \frac{a_{n}}{\lambda_n^3} \right)
\end{eqnarray*}
for large $n$. Hence $\beta_{n}\left( t_{n}(\min(x,y)),t_{n}(\max(x,y)) \right)>0$
and $\tilde{\beta}_{n}\left( t_{n}(\min(x,y)),t_{n}(\max(x,y)) \right)<0$ for large $n$.
Using \eqref{eq3.19}, for large $n$ we have
\begin{eqnarray}\label{eq4.11}
& & \frac{n}{2\pi}\int_{\min\left( \beta_{n}(t_{n}(\min(x,y)), t_{n}(\max(x,y))) , \frac{\pi}{2} \right)}^{\frac{\pi}{2}} \left( 1-F(\lx{t_{n}(\min(x,y))}/\cos \alpha) \right)d\alpha \nonumber\\
&=& \frac{n(2a_{n}|\min(x,y)|)^{\frac{1}{2}}}{4\pi}\int_{\min\left(  \frac{\sin^2\beta_{n}\left( t_{n}(\min(x,y)), t_{n}(\max(x,y)) \right)}{2a_{n}|\min(x,y)|},1-\frac{a_{n}|\min(x,y)|}{2}  \right)}^{1-\frac{a_{n}|\min(x,y)|}{2}}
s^{-\frac{1}{2}}(1-2a_{n}|\min(x,y)|s)^{-\frac{1}{2}}\nonumber \\
& & \times \left( 1-F\left( t_{n}(\min(x,y))(1-2a_{n}|\min(x,y)|s)^{-\frac{1}{2}} \right) \right)
ds  \nonumber \\
&=& 0
\end{eqnarray}
and
\begin{eqnarray}\label{eq4.12}
& & \frac{n}{2\pi}\int_{\max\left( \tilde{\beta}_{n}(t_{n}(\min(x,y)), t_{n}(\max(x,y))) , -\frac{\pi}{2} \right)}^{\frac{\pi}{2}} \left( 1-F(\lx{t_{n}(\max(x,y))}/\cos \alpha) \right)d\alpha \nonumber\\
&=& \frac{n}{2\pi} \int_{0}^{\frac{\pi}{2}} \left( 1-F(\lx{t_{n}(\max(x,y))}/\cos \alpha) \right)d\alpha\nonumber \\
& &
+ \frac{n}{2\pi}  \int_{0}^{\min\left(-\tilde{\beta}_{n}(t_{n}(\min(x,y)), t_{n}(\max(x,y))) ,\frac{\pi}{2}\right)} \left( 1-F(\lx{t_{n}(\max(x,y))}/\cos \alpha) \right)d\alpha   \nonumber \\
&=& |\max(x,y)|^{\frac{1}{2}+\alpha} \left\{
1+\frac{a_{n}(|\max(x,y)|-1)}{c_{\alpha}}\int_{0}^{1}(1-s)^{\alpha-1}s^{\frac{1}{2}}\left( \alpha+1-s-\frac{3}{2}\alpha s \right) ds \right. \nonumber\\
& & \lx{+} A(a_{n}^{-1})\left(\frac{|\max(x,y)|^{-\tau}-1}{c_{\alpha}} \int_{0}^{1}(1-s)^{\alpha}\frac{(1-s)^{-\tau}-1}{\tau}s^{-\frac{1}{2}}ds  \right.
\nonumber\\
&&\left.\left.+ \frac{|\max(x,y)|^{-\tau}-1}{\tau}\right)
+o\left( a_{n}+A(a_{n}^{-1}) \right)  \right\}.
\end{eqnarray}

Combining \eqref{eq4.11} and \eqref{eq4.12}, we have
\begin{eqnarray*}
& & n\Big( 1-\P(\xi_{n1} \leq t_{n}(\min(x,y)), \eta_{n1} \leq t_{n}(\max(x,y))) \Big)- |\min(x,y)|^{\alpha+\frac{1}{2}} \\
&=& |\min(x,y)|^{\frac{1}{2}+\alpha} \left\{
\frac{a_{n}(|\min(x,y)|-1)}{c_{\alpha}}\int_{0}^{1}(1-s)^{\alpha-1}s^{\frac{1}{2}}\left( \alpha+1-s-\frac{3}{2}\alpha s \right) ds \right. \nonumber\\
& &
+A(a_{n}^{-1})\left(\frac{|\min(x,y)|^{-\tau}-1}{c_{\alpha}} \int_{0}^{1}(1-s)^{\alpha}\frac{(1-s)^{-\tau}-1}{\tau}s^{-\frac{1}{2}}ds\right.\nonumber\\
&&\left.\left.+\frac{|\min(x,y)|^{-\tau}-1}{\tau} \right)+o\left( a_{n}+A(a_{n}^{-1}) \right)  \right\}
\end{eqnarray*}
for large $n$, which implies that \eqref{eq2.7} also holds. The proof is complete.
\qed

\noindent
{\bf Proof of Theorem \ref{th7}.}
(i) By $\lim_{n\to \infty} \frac{1-\rho_{n}}{a_{n}}=\infty$ and $\lim_{n\to \infty} \frac{(1-\rho_{n})^3}{a_{n}}=0$,
we have
\begin{eqnarray*}
\frac{\frac{t_{n}(y)}{t_{n}(x)}-\rho_{n}}{\sqrt{a_{n}}\sin\tau_{n}}=\lambda_{n}+\frac{y-x}{2\lambda_{n}}
+\frac{a_{n}\lambda_{n}^3}{2}+o(a_{n}\lambda_{n}^{3})
\end{eqnarray*}
for large $n$, which implies that
\begin{eqnarray}\label{eq4.13}
\frac{\sin\beta_{n}(t_{n}(x),t_{n}(y))}{\sqrt{a_{n}}}=\lambda_{n} +\frac{y-x}{2\lambda_{n}}+o(a_{n}\lambda_{n}^3)
\end{eqnarray}
and
\begin{eqnarray}\label{eq4.14}
\frac{\sin\tilde{\beta}_{n}(t_{n}(x),t_{n}(y))}{\sqrt{a_{n}}}= \lambda_{n}+\frac{x-y}{2\lambda_{n}} +o(a_{n}\lambda_{n}^{3}).
\end{eqnarray}
Here, \eqref{eq4.13} and \eqref{eq4.14} show that $\beta_{n}(t_{n}(x),t_{n}(y))>0$ and
$\tilde{\beta}_{n}(t_{n}(x),t_{n}(y))>0$ for large $n$. Then combining with \eqref{eq3.19},
\begin{eqnarray*}
& & \lx{-}n\P(X\leq t_{n}(x),Y\leq t_{n}(y))+|x|^{\alpha+\frac{1}{2}}+|y|^{\alpha+\frac{1}{2}} \\
&=& -n(1-G(t_{n}(x)))+|x|^{\alpha+\frac{1}{2}} -n(1-G(t_{n}(y)))+|y|^{\alpha+\frac{1}{2}} \\
&=& -|x|^{\frac{1}{2}+\alpha} \left(
\frac{a_{n}(|x|-1)}{c_{\alpha}}\int_{0}^{1}(1-s)^{\alpha-1}s^{\frac{1}{2}}\left( \alpha+1-s-\frac{3}{2}\alpha s \right) ds \right. \nonumber\\
& &  \left.
+A(a_{n}^{-1})\left( \frac{|x|^{-\tau}-1}{\tau} +\frac{|x|^{-\tau}-1}{c_{\alpha}} \int_{0}^{1}(1-s)^{\alpha}\frac{(1-s)^{-\tau}-1}{\tau}s^{-\frac{1}{2}}ds  \right) \right) \\
& & -|y|^{\frac{1}{2}+\alpha} \left(
\frac{a_{n}(|y|-1)}{c_{\alpha}}\int_{0}^{1}(1-s)^{\alpha-1}s^{\frac{1}{2}}\left( \alpha+1-s-\frac{3}{2}\alpha s \right) ds \right. \nonumber\\
& &  \left.
+A(a_{n}^{-1})\left( \frac{|y|^{-\tau}-1}{\tau} +\frac{|y|^{-\tau}-1}{c_{\alpha}} \int_{0}^{1}(1-s)^{\alpha}\frac{(1-s)^{-\tau}-1}{\tau}s^{-\frac{1}{2}}ds  \right) \right)
+o\left( a_{n}+A(a_{n}^{-1}) \right)
\end{eqnarray*}
for large $n$. Hence we can derive \eqref{eq2.3} with $\lambda=\infty$.

(ii) For $\lim_{n\to \infty} \rho_{n}=d\in [-1,1) $, we have
\begin{eqnarray*}
\beta_{n}(t_{n}(x),t_{n}(y))\to \arctan\frac{1-d}{\sin(\arccos d)}=\lx{\frac{1}{2}\arccos d}>0
\end{eqnarray*}
and
\begin{eqnarray*}
\tilde{\beta}_{n}(t_{n}(x),t_{n}(y))\to \arccos d - \arctan\frac{1-d}{\sin(\arccos d)}=\lx{\frac{1}{2}\arccos d}>0
\end{eqnarray*}
as $n\to \infty$.
Using the same arguments as the case of (i), the desired result is derived, which complete
the proof.
\qed

\vspace{1cm}

\noindent {\bf Acknowledgements}~~
This work was supported by the National Natural Science Foundation of China (grant No. 11171275 and No. 11501113),
the Key Project of Fujian Education Committee (grant No. JA15045), the Hujiang Fund for Ph.D. Scientific
Research Startup (grant No. BSQD201608).

\end{document}